\documentclass[11pt]{article}
\usepackage{amsmath}
\usepackage{mathrsfs}
\usepackage{amsfonts}
\usepackage{amssymb}

\newtheorem{lemma}{Lemma}[section]
\newtheorem{proposition}{Proposition}[section]
\newtheorem{theorem}{Theorem}[section]

\newtheorem{remark}{Remark}[section]
\newtheorem{definition}{Definition}[section]

\let\Section=\section
\def\section{\setcounter{equation}{0}\Section}
\sf

\begin{document}
\title{\hspace{-0.2cm}{Schr\"{o}dinger-Poisson system with steep potential
 well}\thanks{This work was supported by
NSFC and CAS-KJCX3-SYW-S03.  \hfil\break\indent
$^\dag$Corresponding author: hszhou@wipm.ac.cn.\hfil\break\indent
{\it 2000 Mathematical Subject Classification}: 35J60.
 {\it Key words}: Schr\" odinger-Poisson, priori
estimate, potential well, ground state, exponential decay. }}
 \author{\small Yongsheng Jiang and Huan-Song Zhou$^\dag$\\
\small Wuhan Institute of Physics and Mathematics\\ \small Chinese
Academy of Sciences, P.O.Box 71010, Wuhan 430071, China
 }
 \date{}
\maketitle
 {\sc Abstract}: We study the following
Schr\"{o}dinger-Poisson system
\begin{equation}\nonumber
(P_\lambda)\left\{\begin{array}{ll}
 -\Delta u + (1+\mu g(x))u+\lambda \phi (x) u =|u|^{p-1}u,\ x\in \mathbb{R}^3, \\
 -\Delta\phi = u^2, \lim\limits_{|x|\rightarrow +\infty}\phi(x)=0,
 \end{array}\right.
\end{equation}
where  $\lambda,\mu$ are positive parameters, $p\in(1,5)$,
$g(x)\in L^\infty(\mathbb{R}^3)$ is nonnegative and $g(x)\equiv0$
on a bounded domain in $\mathbb{R}^3$, $g(x)\not\equiv g(|x|)$.
For $\mu=0$, $\rm(P_\lambda)$ was studied by Ruiz [J. Funct.
Anal., 237(2006), 655-674]. However, if $\mu\neq0$ and $g(x)$ is
not radially symmetric, it is unknown whether $\rm(P_\lambda)$ has
a nontrivial solution for $p\in(1,2)$. In this paper, combining
domain approximation and priori estimate we establish the
boundedness and compactness of a (PS) sequence, then we prove for
$\mu>0$ large that $\rm(P_\lambda)$ with $p\in(1,2)$ has a ground
state if $\lambda>0$ small and that $\rm(P_\lambda)$ for
$p\in[3,5)$ has a nontrivial solution for all $\lambda>0$.
Moreover, some behaviors of the solutions of $\rm(P_\lambda)$ as
$\lambda\rightarrow0$, $\mu\rightarrow+\infty$ and
$|x|\rightarrow+\infty$ are also discussed.

\section{Introduction}
In this paper, we are concerned with the existence of nontrivial
solutions of the following Schr\"{o}dinger-Poisson system
\begin{equation}\label{eq:1.1}
\left\{\begin{array}{ll}
 -\Delta u + V_\mu(x)u+\lambda\phi (x) u =|u|^{p-1}u,\,\,\, x\in \mathbb{R}^3, \\
 -\Delta\phi = u^2,\,\lim\limits_{|x|\rightarrow
 +\infty}\phi(x)=0,\
 \  x\in \mathbb{R}^3,
 \end{array}\right.
\end{equation}
where $V_\mu(x)=1+\mu g(x)$, $\lambda$ and $\mu$ are positive
parameters, $p\in(1,5)$ and $g(x)$ satisfies
\begin{description}
 \item $\rm (G_1)$  $g(x)\geqslant0$, $g(x)\in L^{\infty}(\mathbb{R}^3)$.
\item $\rm (G_2)$ $\Omega_0=\{x\in\mathbb{R}^3:g(x)=0\}$ is
bounded and has nonempty interior. \item $\rm (G_3)$
$\liminf\limits_{|x|\rightarrow\infty}g(x)=1$. (Here 1 can be
replaced by any positive number)
\end{description}
Problem (\ref{eq:1.1}) arisen in quantum mechanics, which is
related to the study of nonlinear Schr\"odinger equation for a
particle in an electromagnetic field or the Hatree-Fock equation,
etc., see \cite{VBenciDFortucto-TopMNA, VBenciDFortucto-RevMP,
D'AprileTMugnaiD-PRSESM, d'AveniaP-ANonSt, RuizD-JFA} and the
references therein. If $\mu=0$, the existence of solutions of
problem (\ref{eq:1.1}) has been discussed for different ranges of
$p$, for examples, \cite{GMCoclite-CAA}
\cite{D'AprileTMugnaiD-PRSESM} \cite{d'AveniaP-ANonSt} for
$p\in[3,5)$, \cite{AzzolliniAandPomponioA-Jmaa} for $p\in(2,5)$,
\cite{KikuchiHiroaki-NonlAnal} for $p\in[2,3)$ and
\cite{AmbrosettiAntonioRuizDavid-CCM}
\cite{RuizD-JFA}\cite{AzzolliniAveniaPomponio-PreP} for
$p\in(1,5)$ or general nonlinearity, etc. For some $\lambda_0>0$,
it was proved in \cite{KikuchiHiroaki-NonlAnal} that
(\ref{eq:1.1}) with $\mu=0$ and $p\in(1,2]$ has no nontrivial
solution if $\lambda\geqslant\lambda_0$. Ruiz in \cite{RuizD-JFA}
proved that $\lambda_0=\frac{1}{4}$ and
 that (\ref{eq:1.1}) has solution if $\lambda>0$ small enough. Some recent results in this direction was summarized in
\cite{AmbrosettiA-MJM}. However, if $\mu\neq0$, there are only a
few results on (\ref{eq:1.1}) for $p\in(2,5)$ and some special
potential functions, such as
\cite{AzzolliniAandPomponioA-Jmaa,ZhaoLeigaZhaoFukun-Jmaa,WangZhengpingZhouHuan-Song-DisCDS,YangShenDing-NA}.
In \cite{AzzolliniAandPomponioA-Jmaa}, a ground state of
(\ref{eq:1.1}) with $p\in(3,5)$ was obtained if $V_\mu(x)$
satisfies
\begin{description}
 \item $\rm (V_1)$  $V_\infty=\lim\limits_{|y|\rightarrow+\infty}V_\mu(y)\geqslant
 V_\mu(x)$
 a.e. in $x\in\mathbb{R}^3$,
 and the strict inequality holds on a positive measure set.
\end{description}
 It was also
 proved in \cite{ZhaoLeigaZhaoFukun-Jmaa} that
(\ref{eq:1.1}) still has a ground state for $p\in(2,3]$ if
$\rm(V_1)$ holds and $V_\mu(x)$ is weakly differentiable function
such that
\begin{description}
 \item $\rm (V_2)$ $(\nabla V_\mu(x),x)\in L^\infty(\mathbb{R}^3)\cup L^{3/2}(\mathbb{R}^3)$ and $2V_\mu(x)+(\nabla V_\mu(x),x) \geqslant0$ a.e.
 $x\in\mathbb{R}^3$,
 where $(\cdot,\cdot)$ is the usual inner product in $\mathbb{R}^3$.
\end{description}
If $-\Delta u$ is replaced by $-\varepsilon \Delta u$ in
(\ref{eq:1.1}), the semiclassical states for this kind of problem
was discussed very recently in \cite{RuizVaira-reprint,
YangShenDing-NA}  for special $V_\mu(x)$, and in
\cite{D'AprileTWei-SIAM} for $\mu=0$, etc.
 Problem (\ref{eq:1.1}) with asymptotically linear nonlinearity was
considered in \cite{WangZhengpingZhouHuan-Song-DisCDS}. To the
authors' knowledge, it seems still open if there exists a solution
to problem (\ref{eq:1.1}) with $\mu\neq0$ and $p\in(1,2)$. It is
known that $\rm (V_1)$ is crucial in using concentration
compactness principle if $V_\mu(x)$ is not radially symmetric,
$\rm(V_2)$ is important in getting the boundedenss of a (PS)
sequence. In this paper, we assume neither the conditions of
$\rm(V_1)$ and $\rm(V_2)$ nor that the potential $V_\mu(x)$ is
radially symmetric. These lead to several difficulties in using
variational method to get solutions of (\ref{eq:1.1}). The first
is to prove that a (PS) sequence is bounded in
$H^1(\mathbb{R}^3)$, which is known to be hard for $p\in(1,2)$
because of the presence of the so called nonlocal term
``$\phi(x)u$"  in (\ref{eq:1.1}). If $\mu=0$, we can overcome this
difficulty by using the decay property of functions in
$H_r^1(\mathbb{R}^3)$ \cite{RuizD-JFA}, but this trick does not
work when $V_\mu(x)$ is not radially symmetric. Once we get a
bounded (PS) sequence, then we have to find a new way to prove
that the (PS) sequence converges to a solution because the
concentration compactness principle cannot be used in our case due
to the lack of conditions $\rm(V_1)$ etc. In this paper, we first
consider the problem (\ref{eq:1.1}) on
$B_k=\{x\in\mathbb{R}^3:|x|<k\}$ for each $k\in\mathbb{N}$, and we
can easily get a solution $u_k$ of (\ref{eq:1.1}) in $H_0^1(B_k)$
since the Sobolev embedding $H_0^1(B_k)\hookrightarrow L^q(B_k)$
for $q\in(1,6)$ is compact. Then, we establish a uniform priori
estimate of $L^\infty$ norm of $u_k$. Finally, by adapting some
techniques used in \cite{WangZhengpingZhouHuan-Song-DisCDS}, we
prove that $\{u_k\}$ converges to a nontrivial solution of
(\ref{eq:1.1}). Moreover, a ground state of (\ref{eq:1.1}) for
$p\in(1,2)$ is also proved in Section 6. However, our methods seem
not useful for $p\in[2,3)$.

{\bf Notations:} Throughout this paper, for $k\in\mathbb{N}$, we
denote $B_k=\{x\in\mathbb{R}^3:|x|<k\}$ and define
\begin{equation}\label{eq:1.2.0}
C^1_B(\mathbb{R}^3)=\{u\in C^1(\mathbb{R}^3): u \in
L^\infty(\mathbb{R}^3) \mbox{ and } |\nabla u|\in
L^\infty(\mathbb{R}^3)\},
\end{equation}
\begin{equation}\label{eq:1.2}
\mathscr{D}_V=\{u\in
\mathscr{D}^{1,2}(\mathbb{R}^3):\int_{\mathbb{R}^3}V_\mu(x)u^2dx<+\infty\},
\end{equation}
\begin{equation}\label{eq:1.3}
\mathscr{D}_k=\{u\in
\mathscr{D}_0^{1,2}(B_k):\int_{B_k}V_\mu(x)u^2dx<+\infty\}.
\end{equation}
 Clearly, $\mathscr{D}_V$ and
$\mathscr{D}_k$ are Hilbert spaces, their scalar products are given
by
\begin{equation}\label{eq:1.4}
\langle u,v\rangle_{\mathscr{D}_V}=\int_{\mathbb{R}^3}\nabla u\nabla
v+V_\mu(x)u vdx,\ \text{ for any } u,v\in \mathscr{D}_V,
\end{equation}
\begin{equation}\label{eq:1.5}
\langle u,v\rangle_k=\int_{B_k}\nabla u\nabla v+V_\mu(x)u vdx, \
\text{ for any } u,v\in \mathscr{D}_k.
\end{equation}
 The norms of $\mathscr{D}_V$ and $\mathscr{D}_k$ are introduced by
\begin{equation}\label{eq:1.7}
\|u\|^2_{\mathscr{D}_V}=\langle u,u\rangle_{\mathscr{D}_V},\ \
\|u\|^2_k=\langle u,u\rangle_k.
\end{equation}
Clearly, $\|\cdot\|_{\mathscr{D}_V}$ is an equivalent norm of
$H^1(\mathbb{R}^3)$. For $p\in[1,+\infty]$, we denote the usual norm
of $L^p(\Omega)$ by $|\cdot|_{L^p(\Omega)}$, and simply by
$|\cdot|_p$ if $\Omega=\mathbb{R}^3$.

Moreover, for any $u\in \mathscr{D}_k$, the extension of $u$ on
$\mathbb{R}^3$ is defined by
\begin{equation}\label{eq:1.8}
\tilde{u}(x)=u(x) \text{ if } x\in B_k,\ \ \ \tilde{u}(x)=0 \text{
if } x\in B^c_k.
\end{equation}
For $u\in \mathscr{D}_V$, let $\phi_u$ be the unique solution of
$-\Delta \phi=u^2$ in $\mathscr{D}^{1,2}(\mathbb{R}^3)$ [Lemma
\ref{Le:2.1} below], define
\begin{equation}\label{eq:1.9}
I(u)=\frac{1}{2}\int_{\mathbb{R}^3}|\nabla u|^2+V(x)u^2
dx+\frac{\lambda}{4}\int_{\mathbb{R}^3}\phi_u(x)u^2dx-\frac{1}{p+1}\int_{\mathbb{R}^3}|u|^{p+1}(x)dx.
\end{equation}
 By
\cite[Proposition 2.1]{D'AprileTMugnaiD-PRSESM}, $I\in
C^1(\mathscr{D}_V,\mathbb{R})$ and for any
$\varphi\in\mathscr{D}_V$ we have
\begin{equation}\label{eq:1.10}
I'(u)\varphi=\int_{\mathbb{R}^3}\nabla u\nabla \varphi+V(x)u\varphi
dx+\lambda\int_{\mathbb{R}^3}\phi_u u\varphi
dx-\int_{\mathbb{R}^3}|u|^{p-1}u\varphi dx.
\end{equation}
Moreover, if $u\in\mathscr{D}_V$ and $I'(u)\varphi=0$ for all
$\varphi\in \mathscr{D}_V$, Lemma 2.4 of
\cite{D'AprileTMugnaiD-PRSESM} showed that $u$ and $\phi_u$
(simply by $\phi$, sometimes) satisfy (\ref{eq:1.1}) in the weak
sense.
\begin{definition}\label{d:1.1}
A pairs $(u,\phi)\in\mathscr{D}_V\times
\mathscr{D}^{1,2}(\mathbb{R}^3)$ is said to be a (weak) solution
of (\ref{eq:1.1}) if
\begin{equation}\label{eq:1.10.0}
I'(u)\varphi=0\ \text{ for all } \varphi\in\mathscr{D}_V.
\end{equation}
\end{definition}
For the sake of simplicity, in many cases we just say
$u\in\mathscr{D}_V$, instead of $(u,\phi_u)\in\mathscr{D}_V\times
\mathscr{D}^{1,2}(\mathbb{R}^3)$, is a weak solution of
(\ref{eq:1.1}).\\ We end this section by giving our main results.
\begin{theorem}\label{th:1.1} Let $p\in(1,2)$ and
$\rm(G_1)$ to $\rm(G_3)$ hold. Then there exists
$\lambda_\ast\in(0,+\infty)$ and $C_{p,\lambda}=2^{\frac{2}{p-2}}
\left(2-p\right)^{\frac{2p-1}{p-1}}\left[p(p-1)\lambda^{-1}\right]^{\frac{p}{2-p}}$
such that problem (\ref{eq:1.1}) has one nontrivial solution
$(u,\phi)\in \mathscr{D}_V\times \mathscr{D}^{1,2}(\mathbb{R}^3)$
for any $\lambda\in(0,\lambda_\ast)$ and $\mu>\mu_1=
C_{p,\lambda}^{p-1}-1$. Moreover,
\begin{equation}\nonumber
I(u)\in[\alpha, c_\lambda],\ \hfill
\|u\|_{\mathscr{D}_V}+|\nabla\phi|_2\leqslant M_\mu:=
M(p,\lambda,\mu), \mbox{ and }
\end{equation}
\begin{equation}\nonumber
 |u(x)|\leqslant C_{p,\lambda}, \hfill \
0<\phi(x)\leqslant C_\mu:=C(p,\lambda,\mu) \ \ a.e.\ \text{in} \
x\in \mathbb{R}^3,
\end{equation}
where $\alpha, c_\lambda>0$ are independent of $\mu$, and $M_\mu$ is
decreasing in $\mu$.
\end{theorem}
\begin{remark}\label{r:1.1} We claim that
$\lambda_\ast\leqslant
c(p)=\frac{1}{4}(p-1)^2(2-p)^{2\frac{2-p}{p-1}}$. In fact, similar
to the proof of \cite[Theorem 4.1]{RuizD-JFA}, we see that
(\ref{eq:1.1}) has no any nontrivial solution in
$\mathscr{D}_V\times\mathscr{D}^{1,2}(\mathbb{R}^3)$ if
$\lambda\geqslant c(p)$. Moreover, it is not difficult to check
that $C_{p,\lambda}>1$ if $\lambda\in(0,c(p))$. Hence $\mu_1>0$ in
Theorem \ref{th:1.1}.
\end{remark}
By Theorem \ref{th:1.1} and Remark \ref{r:1.1}, we know that the
existence of solutions of (\ref{eq:1.1}) for $p\in(1,2)$ depends
heavily on the parameter $\lambda$. However, if $p\in[3,5)$, the
situation is quite different and we have the following theorem.

\begin{theorem}\label{th:1.2} Let $p\in[3,5)$ and $\rm(G_1)$ to
$\rm(G_3)$ hold. Then, for any $\lambda>0$, there exist positive
constants $M_0:=M_0(\lambda)$, $M_1:=M_1(\lambda)$ and
$M:=M(\lambda)$, which are independent of $\mu$ and nondecreasing in
$\lambda>0$, such that problem (\ref{eq:1.1}) has one nontrivial
solution $(u,\phi)\in \mathscr{D}_V\times
\mathscr{D}^{1,2}(\mathbb{R}^3)$ if
$\mu>\mu_2=\max\{0,M_0(\lambda)^{p-1}-1\}$, and $(u,\phi)$ satisfies
\begin{equation}\nonumber
|u(x)|\leqslant M_0\text{ and } \  0<\phi(x)\leqslant M_1\
a.e.\text{ in } x\in \mathbb{R}^3,\
\|u\|_{\mathscr{D}_V}+|\nabla\phi|_2\leqslant M.
\end{equation}
 Moreover, there
exist $\alpha>0$ and $c_\lambda>0$, independent of $\mu$, such
that $I_{\lambda}(u)\in[\alpha, c_\lambda]$ where $\alpha$ is also
independent of $\lambda>0$.
\end{theorem}

The following three theorems describe how the solutions of
(\ref{eq:1.1}) behave when $\lambda\rightarrow0$,
$\mu\rightarrow+\infty$ and $|x|\rightarrow+\infty$, respectively.
\begin{theorem}\label{th:1.3} Under the assumptions of
Theorem \ref{th:1.2}. For each $\lambda>0$ and $\mu>\mu_2$, let
$u_\lambda$ denote the solution of (\ref{eq:1.1}) obtained by
Theorem \ref{th:1.2}. Then there exists a solution $u_0\in
\mathscr{D}_V\setminus \{0\}$ for (\ref{eq:1.1}) with $\lambda=0$
such that, along a subsequence,
\begin{equation}\nonumber
u_\lambda\rightarrow u_0\ \ \text{ strongly in } \mathscr{D}_V \ \
\text{as }\lambda\rightarrow0.
\end{equation}
\end{theorem}

\begin{theorem}\label{th:1.4} For
each $\mu>0$ large, let $u_\mu$ be a solution of (\ref{eq:1.1})
obtained by Theorem \ref{th:1.1} or \ref{th:1.2}. Then, there is
$\bar{u}\in H^1(\mathbb{R}^3)$ with $\bar{u}(x)=0$ a.e. in
$x\in\mathbb{R}^3\setminus\Omega_0$ and $\bar{u}(x)\not\equiv 0$
in $\Omega_0$ such that, passing to a subsequence,
\begin{equation}\nonumber
u_\mu\rightarrow \bar{u}\ \ \text{in }H^1(\mathbb{R}^3) \text{ as
}\mu\rightarrow+\infty.
\end{equation}
Moreover, if $\partial\Omega_0$ is Lipschitz continuous, then
$\bar{u}\in H^1_0(\Omega_0)$ and it is a weak solution of the
following problem
\begin{equation}\label{eq:1.11}
\left\{\begin{array}{ll}
 -\Delta u +u+\frac{\lambda}{4\pi}(u^2\ast\frac{1}{|x|}) u =|u|^{p-1}u,\,\,\, x\in \Omega_0,\\
 u(x)=0,\,\hfill x\in \partial \Omega_0.
 \end{array}\right.
\end{equation}
\end{theorem}

\begin{theorem}\label{th:1.5} For each $\lambda\in(0,\lambda_\ast)$, let $u_\mu$ be the
solutions of (\ref{eq:1.1}) obtained by Theorem \ref{th:1.1} or
\ref{th:1.2}. If
$\mu>\mu_0=\max\{3C_{p,\lambda}^{p-1}-3,3M_0^{p-1}-3\}$,
 then there exist $A>0$ and $R_0>0$ which are independent of
 $\mu$ such that
\begin{equation}\nonumber
u_\mu(x)\leqslant A
|x|^{-\frac{1}{2}}e^{-\frac{\sqrt{\mu}}{2}(|x|-R_0)}\ \text{ for }
|x|>R_0 \text{ and }\mu>\mu_0.
\end{equation}
\end{theorem}
\begin{theorem}\label{th:6.1} For $\lambda_\ast$ and $\mu_1$ given by Theorem \ref{th:1.1}, let $p\in(1,2)$ and
$\rm(G_1)$ to $\rm(G_3)$ hold. Then problem (\ref{eq:1.1}) has a
ground state $(u,\phi)\in \mathscr{D}_V\times
\mathscr{D}^{1,2}(\mathbb{R}^3)$ for each
$\lambda\in(0,\lambda_\ast)$ and $\mu>\mu_1$.
\end{theorem}

\section{ Nontrivial solution for (\ref{eq:1.1}) on $B_k$}
For any $k\in\mathbb{N}$, we consider the following problem
\begin{equation}\label{eq:2.1}
\left\{\begin{array}{ll}
 -\Delta u + V_\mu(x)u+\lambda\phi (x) u =|u|^{p-1}u,\,\,\, x\in B_k,\\
 -\Delta\phi = u^2,\ \hfill x\in  B_k,\\
 u(x)=\phi(x)=0,\,\hfill x\in \partial B_k.
 \end{array}\right.
\end{equation}
For $u\in \mathscr{D}_k$ and
$\phi:=\phi_u\in\mathscr{D}^{1,2}_0(B_k) $, define
\begin{equation}\label{eq:2.2}
I_{k}(u):=I_{\lambda,k}(u)\triangleq\frac{1}{2}\int_{B_k}|\nabla
u|^2+V_\mu(x)u^2
dx+\frac{\lambda}{4}\int_{B_k}\phi_u(x)u^2dx-\frac{1}{p+1}\int_{B_k}|u|^{p+1}dx,
\end{equation}
and $I_{k}\in C^1(\mathscr{D}_k,\mathbb{R})$,
$(u,\phi_u)\in\mathscr{D}_k\times\mathscr{D}^{1,2}_0(B_k)$ is a weak
solution of (\ref{eq:2.1}) if and only if $u$ is a nonzero critical
point of $I_k(u)$. The main aim of this section is to prove that
\begin{theorem}\label{th:2.1} Let $p\in(1,2)\cup[3,5)$ and $\rm(G_1)$ to
$\rm(G_3)$ hold. Then there exist $\lambda_\ast>0$ and
$k_0\in\mathbb{N}$ such that, for $k>k_0$, problem (\ref{eq:2.1})
has at least one nontrivial solution $(u_k,\phi_k)\in
\mathscr{D}_k\times \mathscr{D}_0^{1,2}(B_k)$ with
\begin{equation}\nonumber
I_{k}(u_k)=c_{\lambda,k}\in[\alpha, c_\lambda]\ \text{for each
}\lambda\in(0,\lambda_\ast),
\end{equation}
where $\lambda_\ast=+\infty$ if $p\in[3,5)$ and
$\lambda_\ast\in(0,+\infty)$ if $p\in(1,2)$, $\alpha$ and
$c_\lambda$ are positive constants independent of $k$ and $\mu$ (
$\alpha$ is also independent of $\lambda$). Moreover, for $k>k_0$,
$c_{\lambda_1,k}\leqslant c_{\lambda_2,k}$ if
$\lambda_1\leqslant\lambda_2$.
\end{theorem}

Before proving Theorem \ref{th:2.1}, we recall some results for the
following Poisson equation on any smooth domain
$\Omega\subset\mathbb{R}^3$, $\Omega$ may be unbounded,
\begin{equation}\label{eq:2.3}
\left\{\begin{array}{ll}
 -\Delta \phi(x) =u^2,\,\,\, \ x\in \Omega, \\
 \phi\in \mathscr{D}_0^{1,2}(\Omega).
 \end{array}\right.
\end{equation}
The following lemma is a result of Lax-Milgram theorem.
\begin{lemma}\label{Le:2.1}
Let $u\in L^{12/5}(\Omega)$. Then (\ref{eq:2.3}) has a unique
solution $\phi_u\in \mathscr{D}_0^{1,2}(\Omega)$  such that
\begin{equation}\label{eq:2.4}
|\nabla\phi_u|_{L^2(\Omega)}\leqslant
S_0^{-1/2}|u|^2_{L^\frac{12}{5}(\Omega)},\ \ \
\int_{\Omega}\phi_uu^2(x)dx\leqslant
S_0^{-1}|u|^4_{L^\frac{12}{5}(\Omega)},
\end{equation}
where $S_0=\inf\{\int_{\mathbb{R}^3}|\nabla u|^2dx:u\in
\mathscr{D}^{1,2}(\mathbb{R}^3)\ \text{ and }
\int_{\mathbb{R}^3}|u|^6dx=1 \}$ is the Sobolev constant, which is
independent of $\Omega$ and $u$.
\end{lemma}
\begin{lemma}\label{Le:2.2}(\cite[Lemma 2.1]{RuizD-JFA})
For $\{u_n\}\subset L^{12/5}(\Omega)$ with
$u_n\overset{n}{\rightarrow}u$ strongly in $L^{12/5}(\Omega)$, let
$\phi_n$ and $\phi$ be the unique solutions of (\ref{eq:2.3})
corresponding to $u_n$ and $u$, respectively. Then,
\begin{equation}\nonumber
\phi_n\overset{n}{\rightarrow}\phi \ \ \text{strongly in }
\mathscr{D}_0^{1,2}(\Omega) \ \text{ and }\ \int_{\Omega} \phi_n
u_n^2 dx\overset{n}{\rightarrow} \int_{\Omega} \phi u^2dx.
\end{equation}
\end{lemma}
\begin{lemma}\label{Le:2.3}
Let $\{u_k\}\subset L^{12/5}(\mathbb{R}^3)$ be such that
$u_k\overset{k}{\rightarrow}u$ strongly in $L^{12/5}(\mathbb{R}^3)$
for some $u\in L^{12/5}(\mathbb{R}^3)$, and $\phi_k\in
\mathscr{D}_{0}^{1,2}(B_k)$ be the unique solution of (\ref{eq:2.3})
with $\Omega=B_k$ and $u=u_k$. If $\tilde{\phi}_k$ is the extension
of $\phi_k$ on $\mathbb{R}^3$ defined as (\ref{eq:1.8}),
 then
\begin{equation}\nonumber
\tilde{\phi}_k(x)\overset{k}{\rightarrow}\phi(x)=\frac{1}{4\pi}\int_{\mathbb{R}^3}\frac{u^2(y)}{|x-y|}dy\
\text{ strongly in } \mathscr{D}^{1,2}(\mathbb{R}^3),
\end{equation}
where $\phi$ is essentially the unique solution of (\ref{eq:2.3})
with $\Omega=\mathbb{R}^3$.
\end{lemma}
{\it Proof.} For $u\in L^{12/5}(\mathbb{R}^3)$, by Lemma
\ref{Le:2.1} there is a unique $\phi\in
\mathscr{D}^{1,2}(\mathbb{R}^3)$ such that
\begin{equation}\label{eq:2.5}
\int_{\mathbb{R}^3}\nabla\phi\nabla\varphi dx=\int_{\mathbb{R}^3}
u^2\varphi dx, \ \text{ for any } \varphi\in
\mathscr{D}^{1,2}(\mathbb{R}^3).
\end{equation}
Moreover, $\phi(x)=\frac{1}{4\pi}
\int_{\mathbb{R}^3}\frac{u^2(y)}{|x-y|}dy$ by Theorem 9.9 of
\cite{DGilbargNStrudinger}. By Lemma \ref{Le:2.1}, there is $C>0$
independent of $B_k$ and $u_k$ such that
\begin{equation}\nonumber
|\nabla\tilde{\phi}_k|_2=\left\{\int_{B_k}|\nabla\phi_k|^2dx\right\}^{1/2}\leqslant
C|u_k|^2_{L^{12/5}(B_k)}\leqslant C|u_k|^2_{12/5},
\end{equation}
this implies that $\{\tilde{\phi}_k\}$ is bounded in
$\mathscr{D}^{1,2}(\mathbb{R}^3)$, since
$u_k\overset{k}{\rightarrow}u$ strongly in $L^{12/5}(\mathbb{R}^3)$
and $\{u_k\}$ is bounded in $L^{12/5}(\mathbb{R}^3)$. Hence, there
exists $\tilde{\phi}\in \mathscr{D}^{1,2}(\mathbb{R}^3)$ such that $
\tilde{\phi}_k\overset{k}{\rightharpoonup}\tilde{\phi}$ weakly in $
\mathscr{D}^{1,2}(\mathbb{R}^3)$,
 that is,
\begin{equation}\label{eq:2.6}
\int_{\mathbb{R}^3}\nabla\tilde{\phi}_k\nabla\varphi
dx\overset{k}{\rightarrow}
\int_{\mathbb{R}^3}\nabla\tilde{\phi}\nabla\varphi dx, \ \ \text{for
any } \varphi\in \mathscr{D}^{1,2}(\mathbb{R}^3).
\end{equation}
By the definition of $\tilde{\phi}_k$, (\ref{eq:2.6}) yields
\begin{equation}\label{eq:2.7}
\int_{B_k}\nabla\phi_k\nabla\varphi
dx=\int_{\mathbb{R}^3}\nabla\tilde{\phi}_k\nabla\varphi
dx\overset{k}{\rightarrow}
\int_{\mathbb{R}^3}\nabla\tilde{\phi}\nabla\varphi dx,\text{ for any
}\varphi\in C^\infty_0(\mathbb{R}^3).
\end{equation}
 For any $\varphi\in
C^\infty_0(\mathbb{R}^3)$, since
$B_k\overset{k}{\rightarrow}\mathbb{R}^3$ we see that
$\text{supp}\varphi\subset B_k$ for all $k\in\mathbb{N}$ large.
Noting that $\phi_k$ satisfies (\ref{eq:2.3}) with $u=u_k$ and
$\Omega=B_k$, it follows from $u_k\overset{k}{\rightarrow}u$ that
\begin{equation}\label{eq:2.8}
\int_{B_k}\nabla\phi_k\nabla\varphi dx=\int_{B_k}u^2_k\varphi
dx+o(1)=\int_{\mathbb{R}^3} u^2_k\varphi dx+o(1)=\int_{\mathbb{R}^3}
u^2\varphi dx+o(1),
\end{equation}
where and in what follows $o(1)$ denotes the quantity which goes to
0 as $k\rightarrow+\infty$.
 Therefore, combining (\ref{eq:2.7})
(\ref{eq:2.8}) and (\ref{eq:2.5}), we see that
\begin{equation}\nonumber
\int_{\mathbb{R}^3} \nabla\tilde{\phi}\nabla\varphi
dx=\int_{\mathbb{R}^3} u^2\varphi dx=\int_{\mathbb{R}^3}
\nabla\phi\nabla\varphi dx, \text{ for any }\varphi\in
\mathscr{D}^{1,2}(\mathbb{R}^3).
\end{equation}
So, $\tilde{\phi}\equiv\phi$ in $\mathscr{D}^{1,2}(\mathbb{R}^3)$.
Taking $\varphi=\tilde{\phi}_k$ and $\phi$ in (\ref{eq:2.5}) and
(\ref{eq:2.6}), respectively, then,
\begin{equation}\label{eq:2.9}
\int_{\mathbb{R}^3} u^2\tilde{\phi}_k dx=\int_{\mathbb{R}^3}
\nabla\phi\nabla\tilde{\phi}_k
dx\overset{k}{\rightarrow}\int_{\mathbb{R}^3} |\nabla\phi|^2dx.
\end{equation}
On the other hand,
\begin{equation}\nonumber
\int_{\mathbb{R}^3} |\nabla\tilde{\phi}_k|^2 dx=\int_{B_k}
|\nabla\phi_k|^2 dx=\int_{B_k} u_k^2 \phi_k dx=\int_{\mathbb{R}^3}
u_k^2\tilde{\phi}_k dx.
\end{equation}
By $u_k\overset{k}{\rightarrow}u$ strongly in
$L^{12/5}(\mathbb{R}^3)$, it follows from (\ref{eq:2.9}) that
\begin{equation}\nonumber
\int_{\mathbb{R}^3} |\nabla\tilde{\phi}_k|^2
dx\overset{k}{\rightarrow}\int_{\mathbb{R}^3} |\nabla\phi|^2 dx,
\end{equation}
this implies that $\tilde{\phi}_k\overset{k}{\rightarrow}\phi$ in
$\mathscr{D}^{1,2}(\mathbb{R}^3)$. $\Box$\\

Now, we are ready to prove Theorem \ref{eq:2.1}. For using
Mountain Pass Theorem, the following properties for $I_{k}$
defined by (\ref{eq:2.2}) are required.
\begin{lemma}\label{Le:2.4} Let $p\in(1,2)\cup[3,5)$ and $\rm(G_1)$
to $\rm(G_2)$ hold. We have that
\begin{description}
\item $\bf (i)$ There are $\rho>0$ and $\alpha>0$ (independent of
$k$, $\lambda$ and $\mu$) such that
\begin{equation}\nonumber
\inf\limits_{u\in \mathscr{D}_k;
\|u\|_k=\rho}I_{\lambda,k}(u)\geqslant\alpha>0,\text{ for all
}\lambda>0,k\in\mathbb{N}.
\end{equation}
\item $\bf (ii)$ There exist $\lambda_\ast\in(0,+\infty)$
if $p\in(1,2)$ or $\lambda_\ast=+\infty$ if $p\in[3,5)$ and
$e\in\mathscr{D}^{1,2}(\mathbb{R}^3)$ such that $\text{\rm
supp}e\subset B_k$,  $\|e\|_k>\rho$ and
$$
I_{\lambda,k}(e)<0, \text{ for } \lambda\in(0,\lambda_\ast) \text{
and } k\in\mathbb{N} \text{ large},$$ where $e$ is independent of
$\lambda$ only if $p\in(1,2)$.

\item $\bf (iii)$ There is a constant $c_\lambda>0$, independent
of $k$ and $\mu$, such that $$c_{\lambda,k}:
=\inf\limits_{\gamma\in\Gamma_\lambda}\max\limits_{t\in[0,1]}
I_{\lambda,k}(\gamma(t))\leqslant c_\lambda<+\infty,$$ where
$\Gamma_{\lambda}=\{\gamma\in C([0,1], \mathscr{D}_k):\gamma(0)=0,
\|\gamma(1)\|_k>\rho, I_{\lambda,k}(\gamma(1))<0\}$. Moreover, for
any fixed $k\in\mathbb{N}$, $0<\alpha\leqslant
c_{\lambda_1,k}\leqslant c_{\lambda_2,k}$ if
$\lambda_1\leqslant\lambda_2$.
\end{description}
\end{lemma}
{\it Proof. }$\bf(i)$ For $q\in(2,6)$, let
$$S=\inf\{\int_{\mathbb{R}^3}|\nabla
u|^2+u^2dx:u\in H^1(\mathbb{R}^3)\text{ and } |u|_{q}=1\},$$ which
is independent of $k$ and $\mu$. For any $u\in \mathscr{D}_k\subset
H^1_0(B_k)$, define $\tilde{u}$ as (\ref{eq:1.8}) and
$\tilde{u}\in\mathscr{D}_V\subset H^1(\mathbb{R}^3)$. Then, Sobolev
embedding implies that
\begin{equation}\nonumber
\int_{B_k}|u|^{p+1}dx=\int_{\mathbb{R}^3}|\tilde{u}|^{p+1}dx\leqslant
S^{-\frac{p+1}{2}}\|\tilde{u}\|_{H^1(\mathbb{R}^3)}^{p+1}\leqslant
S^{-\frac{p+1}{2}}\|\tilde{u}\|_{\mathscr{D}_V}^{p+1}=S^{-\frac{p+1}{2}}\|u\|_k^{p+1}.
\end{equation}
Since $\lambda>0$ and $\phi_u>0$, it follows from (\ref{eq:2.2})
that
\begin{equation}\nonumber
I_{\lambda,k}(u)\geqslant\frac{1}{2}\|u\|_k^2-\frac{1}{p+1}S^{-\frac{p+1}{2}}\|u\|_k^{p+1}.
\end{equation}
By $p>1$, it is not difficult to see that there exist positive
constants $\rho$, $\alpha$ (independent of $k$, $\lambda$ and $\mu$)
such that $\rm(i)$ holds.\\
\indent$\bf(ii)$ We separate the proof into two cases: $p\in(1,2)$
and $p\in[3,5)$.

 If $p\in(1,2)$, let $w\in
C_0^\infty(\mathbb{R}^3)\setminus\{0\}$ with $\text{supp}w\subset
\Omega_0$. Since $\Omega_0\subset B_k$ for $k\in\mathbb{N}$ large,
there is $k_0\in\mathbb{N}$ such that, for all $k>k_0$ and $t>0$,
\begin{equation}\label{eq:2.10}
I_{0,k}(tw)\triangleq
I_{\lambda,k}(tw)|_{\lambda=0}=\frac{1}{2}\int_{\Omega_0}|\nabla
w|^2+w^2dx-\frac{t^{p+1}}{p+1}\int_{\Omega_0}|w|^{p+1}dx.
\end{equation}
Clearly, there is $t_0>0$ large enough such that $e=t_0w$
satisfying $I_{0,k}(e)<0$ and
$\|e\|_k=t_0\|w\|_{H_0^{1}(\Omega_0)}>\rho$, here $e$ is
independent of $k$, $\lambda$ and $\mu$. Since
$I_{\lambda,k}(e)=I_{0,k}(e)+\frac{\lambda}{4}\int_{\Omega_0}\phi_e
e^2dx$ is continuous in $\lambda\geqslant0$, combining
$I_{0,k}(e)<0$ and Lemma \ref{Le:2.1} we see that there exists
$\lambda_\ast>0$ small (independent of $k$ and $\mu$) such that
$I_{\lambda,k}(e)<0$ for all $\lambda\in(0,\lambda_\ast)$.

 If $p\in[3,5)$, the above proof for $p\in(1,2)$
still works but by that way $\lambda$ has to be small. However,
for $p\in[3,5)$, we can prove that $\rm(ii)$ holds for all
$\lambda>0$, that is $\lambda_\ast=+\infty$. In fact, by $\rm
int\Omega_0\neq\emptyset$ we may assume that
$B_{\varepsilon_0}(x_0)=\{x\in\mathbb{R}^3:|x-x_0|<\varepsilon_0\}\subset\Omega_0$
for some $x_0\in \Omega_0$ and $\varepsilon_0>0$. Taking $w\in
C^\infty_0(\mathbb{R}^3)$ with ${\rm supp}w\subset
B_{\varepsilon_0}(0)$ and letting $w_t(x)=t^2w(t(x-x_0))$, then
${\rm supp} w_t\subset B_{\varepsilon_0}(x_0)$ for $t>1$, and
$B_{\varepsilon_0}(x_0)\subset\Omega_0\subset B_k$ if $k\in
\mathbb{N}$ large enough. Since $g(x)\equiv0$ on $\Omega_0$, for
$t>1$, we have
\begin{align*}
I_{\lambda,k}(w_t)&=\frac{1}{2}\int_{\Omega_0}|\nabla w_t|^2+w_t^2
dx+\frac{\lambda }{4}\int_{\Omega_0}\phi_{w_t}w_t^2dx-\frac{1}{p+1}
\int_{\Omega_0}|w_t|^{p+1}dx\\
 &= \frac{t^3}{2}\int_{\Omega_0}|\nabla w|^2dx+ \frac{t}{2}\int_{\Omega_0}w^2
dx+\frac{\lambda
t^3}{4}\int_{\Omega_0}\phi_{w}w^2dx-\frac{t^{2p-1}}{p+1}
\int_{\Omega_0}|w|^{p+1}dx.
\end{align*}
Since $p>2$, $2p-1>3$, we see that, for each $\lambda>0$, there
exists $t_0:=t_0(\lambda)>1$ large (independent of $k$ and $\mu$)
such that
$$I_{\lambda,k}(w_{t_0})<0 \text{ and } \|w_{t_0}\|_k^2=\frac{t_0^3}{2}\int_{\Omega_0}|\nabla w|^2dx+
\frac{t_0}{2}\int_{\Omega_0}w^2 dx>\rho^2.$$ Hence $\rm(ii)$ is
proved by taking $e:=e_\lambda\triangleq w_{t_0}$ (independent of
$k$ and $\mu$).

$\bf(iii)$ By the result of part $\rm(i)$, it is obvious that
$c_{\lambda,k}\geqslant\alpha>0$.

 For $p\in(1,2)\cup[3,5)$, by part $\rm(ii)$, there
always exists $e\in \mathscr{D}_k$ such that $I_{\lambda,k}(e)<0$.
Let $\gamma_0(t)=te$ for $t\in[0,1]$. Then
$\gamma_0(t)\in\Gamma_\lambda$ and
\begin{equation}\nonumber
c_{\lambda,k}\leqslant\max\limits_{t\in[0,1]}I_{\lambda,k}(\gamma_0(t))=\max\limits_{t\in[0,1]}I_{\lambda,k}(te)\overset{\vartriangle}{=}c_\lambda<+\infty.
\end{equation}
Note that $\text{supp}e\subset\Omega_0$ and $e$ is independent of
$k$ and $\mu$, so $c_\lambda$ is also independent of $k$ and
$\mu$. For each $k>k_0$ and for any $u\in \mathscr{D}_k$,
(\ref{eq:2.2}) implies that $I_{\lambda_1,k}(u)\leqslant
I_{\lambda_2,k}(u)$ if $\lambda_1\leqslant\lambda_2$, then the
definition of $\Gamma_\lambda$ yields that
$\Gamma_{\lambda_2}\subseteq\Gamma_{\lambda_1}$. Therefore
\begin{equation}\nonumber
c_{\lambda_1,k}=\inf\limits_{\gamma\in\Gamma_{\lambda_1}}\max\limits_{t\in[0,1]}
I_{\lambda_1,k}(\gamma(t))\leqslant\inf\limits_{\gamma\in\Gamma_{\lambda_2}}\max\limits_{t\in[0,1]}
I_{\lambda_2,k}(\gamma(t))=c_{\lambda_2,k}.\ \ \Box
\end{equation}

{\bf Proof of Theorem \ref{th:2.1}:} Since $I_k(0)=0$, it follows
from Lemma \ref{Le:2.4} and the mountain pass theorem that, for each
$k\in\mathbb{N}$ large and $\lambda\in(0,\lambda_\ast)$, there
exists $\{u_n\}:=\{u_n^{(k)}\}\subset \mathscr{D}_k$,
$n=1,2,\cdots,$ such that
\begin{equation}\nonumber
I_{k}(u_n)\overset{n}{\rightarrow}c_{\lambda,k}\ \  \text{ and }\ \
I'_{k}(u_n)\overset{n}{\rightarrow}0\ \text{ in }
\mathscr{D}_k^\ast,
\end{equation}
where $\mathscr{D}^\ast_k$ is the dual space of $\mathscr{D}_k$.
Hence, as $n\rightarrow+\infty$,
\begin{equation}\label{eq:2.11}
\frac{1}{2}\int_{B_k}|\nabla
u_n|^2+V_\mu(x)u_n^2dx+\frac{\lambda}{4}\int_{B_k}\phi_nu_n^2dx-\frac{1}{p+1}\int_{B_k}|u_n|^{p+1}dx
=c_{\lambda,k}+o(1),
\end{equation}
\begin{equation}\label{eq:2.12}
 I'_{k}(u_n)u_n =\int_{B_k}|\nabla
u_n|^2+V_\mu(x)u_n^2dx+\lambda\int_{B_k}\phi_nu_n^2dx-\int_{B_k}|u_n|^{p+1}dx,
\end{equation}
\begin{equation}\label{eq:2.13}
\int_{B_k}\nabla u_n\nabla\varphi+V_\mu(x)u_n\varphi
dx+\lambda\int_{B_k}\phi_nu_n\varphi
dx-\int_{B_k}|u_n|^{p-1}u_n\varphi dx=o(1),
\end{equation}
for any $\varphi\in \mathscr{D}_k$, where $\phi_n$ denotes the
unique solution of Poisson equation $-\Delta \phi_n=u_n^2$ in
$\mathscr{D}_k$. For each $k\in\mathbb{N}$ large, if
$\{u_n\}=\{u_n^{(k)}\}$ is bounded in $\mathscr{D}_k$, then there
exists $u_k\in \mathscr{D}_k$ such that, passing to a subsequence,
$$u_n\overset{n}{\rightharpoonup}u_k \text{ weakly in } \mathscr{D}_k \text{ and
strongly in } L^q(B_k) \text{ for } q\in(1,6),$$ since the embedding
$\mathscr{D}_k\subset H_0^1(B_k)\hookrightarrow L^q(B_k)$ is
compact. Hence, Lemma \ref{Le:2.2} follows that
\begin{equation}\nonumber
\int_{B_k}\phi_nu_n^2dx\overset{n}{\rightarrow}\int_{B_k}\phi_{u_k}u_k^2dx,\
\text{ where } -\Delta \phi_{u_k}=u_k^2 \ \text{ in }
\mathscr{D}_0^{1,2}(B_k).
\end{equation}
Based on these facts, it is not difficult to see that
$u_n\overset{n}{\rightarrow}u_k$ strongly in $\mathscr{D}_k$ by
(\ref{eq:2.12}) and (\ref{eq:2.13}). Hence,
$I_k(u_k)=c_{\lambda,k}\in[\alpha,c_\lambda]$ by Lemma \ref{Le:2.4},
and $u_k\not\equiv0$ is a weak solution of (\ref{eq:2.1}).
 Therefore, to prove Theorem \ref{th:2.1} we need only
to show that $\{u_n\}$ is bounded in $\mathscr{D}_k$.\\
 \indent If $p\in
[3,5)$, then $\frac{1}{4}-\frac{1}{p+1}\geqslant0$. Since
$I'_{k}(u_n)\overset{n}{\rightarrow}0$ in $\mathscr{D}_k^\ast$,
subtracted $\frac{1}{4}\times(\ref{eq:2.12})$ from
$(\ref{eq:2.11})$ which gives
\begin{equation}\label{eq:2.14.0}
\|u_n\|_k^2\leqslant 8c_k+2, \ \text{ for } n \text{ large enough, }
\end{equation}
 hence $\{u_n\}$ is bounded in
$\mathscr{D}_k$. However, if $p\in(1,2)$, the properties of $I_k(u)$
change greatly. In this case, we have to use other trick to get the
boundedness of $\{u_n\}$ in $\mathscr{D}_k$. By $-\Delta \phi_n=
u_n^2$ in $\mathscr{D}_0^{1,2}(B_k)$, we have that
\begin{align}\label{eq:2.14}
\sqrt{\lambda}\int_{B_k}|u_n|^3dx&=\sqrt{\lambda}\int_{B_k}\nabla
\phi_n\nabla |u_n|dx \leqslant\frac{1}{2}\int_{B_k}|\nabla
u_n|^2dx+\frac{\lambda}{2}\int_{B_k}|\nabla\phi_n|^2dx \nonumber\\
&=\frac{1}{2}\int_{B_k}|\nabla
u_n|^2dx+\frac{\lambda}{2}\int_{B_k}\phi_n |u_n|^2dx.
\end{align}
Since $I'_k(u_n)\overset{n}{\rightarrow}0$ in $\mathscr{D}_k^\ast$,
we may assume that $\|I'_k(u_n)\|_{\mathscr{D}_k^\ast}\leqslant1$
for $n$ large. Then, by (\ref{eq:2.12}) and (\ref{eq:2.14}) we see
that, for $n\in\mathbb{N}$ large,
\begin{equation}\label{eq:2.15}
\frac{1}{2}\int_{B_k}|\nabla u_n|^2+V_\mu
u_n^2dx+\sqrt{\lambda}\int_{B_k}|u_n|^3dx\leqslant\int_{B_k}|
u_n|^{p+1}dx+\|u_n\|_k.
\end{equation}
If $p\in(1,2)$ and the Young's inequality yields that
\begin{equation}\nonumber
\begin{split}
\int_{B_k}|
u_n|^{p+1}dx&\leqslant\int_{B_k}\left[\frac{(p+1)\varepsilon}{3}|u_n|^3+\frac{2-p}{3}\varepsilon^{-\frac{p+1}{2-p}}\right]dx\\
&=\sqrt{\lambda}\int_{B_k}|u_n|^3dx+\frac{2-p}{3}
\left(\frac{3\sqrt{\lambda}}{p+1}\right)^{-\frac{p+1}{2-p}}|B_k|,\
\text{by taking } \varepsilon=\frac{3\sqrt{\lambda}}{p+1}.
\end{split}
\end{equation}
This and (\ref{eq:2.15}) imply that
\begin{equation}\nonumber
\frac{1}{2}\|u_n\|_k^2\leqslant\|u_n\|_k+\frac{2-p}{3}
\left(\frac{3\sqrt{\lambda}}{p+1}\right)^{-\frac{p+1}{2-p}}|B_k|,\
\text{for } n \text{ large}.
\end{equation}
So, $\{u_n\}=\{u_n^{(k)}\}$ is bounded in $\mathscr{D}_k$ for each
$k\in\mathbb{N}$ large. $\Box$
\section{$L^\infty$-norm priori estimates for solutions of problem (\ref{eq:2.1})}
The aim of this section is to establish priori estimates of
$L^\infty$-norm for solutions of (\ref{eq:1.1}) and
(\ref{eq:2.1}). If $p\in(1,2)$, Lemma \ref{Le:3.1} shows that the
$L^\infty$-norm of solutions of (\ref{eq:1.1}) or (\ref{eq:2.1})
is bounded above by a constant depending only on $\lambda$ and
$p$. If $p\in[3,5)$, similar to (\ref{eq:2.14.0}), we know that
the solution $u_k$ of (\ref{eq:2.1}) is uniformly bounded in
$\mathscr{D}_k$ for all $k\in \mathbb{N}$ and so does in
$L^{6}(B_k)$, then in Lemma \ref{Le:3.3} we establish an estimate
of $|u_k|_{L^\infty(B_k)}$ by using $|u_k|_{L^{6}(B_k)}$. These
priori estimates are important in showing that the solution $u_k$
of (\ref{eq:2.1}) is uniformly bounded with respect to the domain
$B_k$ and in showing that $\tilde{u}_k$, the extension of $u_k$ on
$\mathbb{R}^3$, converges to a nontrivial solution of
(\ref{eq:1.1}).

\begin{lemma}\label{Le:3.1}
For $p\in(1,2)$ and $\Omega=B_k$ or $\Omega=\mathbb{R}^3$, let
$(u,\phi)\in H^1_0(\Omega)\times \mathscr{D}^{1,2}_0(\Omega)$ be a
weak solution of the following problem
\begin{equation}\label{eq:3.1.L}
\left\{\begin{array}{ll}
 -\Delta u + V(x)u+\lambda\phi (x) u =|u|^{p-1}u,\,\,\, x\in \Omega, \\
 -\Delta\phi = u^2,\,\, x\in \Omega,\\
 \end{array}\right.
\end{equation}
 where $\lambda>0$, $V(x)\in L^\infty(\Omega)$ and $V(x)\geqslant1$. Then
\begin{equation}\nonumber
|u(x)|\leqslant c_p\phi(x)\ \text{and } \ |u(x)|\leqslant
C_{p,\lambda},\ \
  \text{a.e. in } x\in \Omega,
\end{equation}
where $c_p=(p-1)(2-p)^{\frac{2-p}{p-1}}$ and
$C_{p,\lambda}=\frac{1}{2}(2-p)\left(\frac{pc_p}{2\lambda}\right)^{\frac{p}{2-p}}$.
\end{lemma}
{\it Proof:} By assumption, $(u,\phi)\in
 H_0^1(\Omega)\times \mathscr{D}^{1,2}_0(\Omega)$  is a
 weak solution of
(\ref{eq:3.1.L}), then, for any $v\in H_0^1(\Omega)$, we have
\begin{equation}\label{eq:3.1}
\int_{\Omega} \nabla u \nabla
vdx+\int_{\Omega}V(x)uvdx+\lambda\int_{\Omega}\phi(x)uvdx
 -\int_{\Omega}|u|^{p-1}u vdx=0,
\end{equation}
\begin{equation}\label{eq:3.2}
\int_{\Omega} \nabla \phi \nabla vdx
 =\int_{\Omega} u^{2} vdx.
\end{equation}
Adding $c_p\int_{\Omega}u^2vdx$ on both sides of (\ref{eq:3.1}),
and using (\ref{eq:3.2}) we get that
\begin{equation}\label{eq:3.2.0}
\begin{split}
\int_{\Omega} \nabla u \nabla
vdx&+\int_{\Omega}[V(x)u+c_pu^2-|u|^{p-1}u]v
dx+\lambda\int_{\Omega}\phi (x)uvdx\\
&=c_p\int_{\Omega} \nabla \phi\nabla vdx,\ \text{ for any } v\in
H_0^1(\Omega).
 \end{split}
\end{equation}
Based on this observation, we prove our lemma by the following two cases.  \\

{\bf Cases 1:} $\Omega=B_k$. We let
\begin{equation}\label{eq:3.2.1}
w_1(x)=(u(x)-c_p\phi(x))^+ \ \text{ and }
\Omega_1=\{x\in\Omega:w_1(x)
>0\}.
\end{equation}
then $w_1 \in H^1_0(\Omega)$ and $u(x)|_{\Omega_1}\geqslant
c_p\phi(x)>0$. Taking $v(x)=w_1(x)$ in (\ref{eq:3.2.0}), and using
$V(x)\geqslant1$, we see that
\begin{equation}\label{eq:3.3.3.1}
\int_{\Omega_1} \nabla u \nabla w_1
dx+\int_{\Omega_1}[u+c_pu^2-|u|^{p-1}u]w_1dx\leqslant
 c_p\int_{\Omega_1} \nabla \phi\nabla w_1 dx.
\end{equation}
However, for all $t\geqslant0$ we have $t+c_pt^2-t^p\geqslant0$ if
$p\in(1,2)$ and $c_p=(p-1)(2-p)^{\frac{2-p}{p-1}}$. Then,
(\ref{eq:3.3.3.1}) implies that $\int_{\Omega_1} \nabla u \nabla
w_1 dx-
 c_p\int_{\Omega_1} \nabla \phi\nabla w_1 dx\leqslant0$, that is,
\begin{equation}\label{eq:3.3.3.2}
 \int_{\Omega_1}\nabla (u-c_p\phi)\nabla w_1dx=
 \int_{\Omega_1} |\nabla w_1|^2dx=
 0.
\end{equation}
Hence, $|\Omega_1|=0$ or $w_1|_{\Omega_1}\equiv\text{constant}$.
By the definition of $\Omega_1$ and $w_1 \equiv 0$ in $\Omega
\setminus \Omega_1$, then $u(x)\leqslant c_p\phi(x)$ a.e. in
$\Omega$. On the other hand, replacing $u$ by $-u$ and repeating
the above procedure, we see that $-u(x)\leqslant c_p\phi(x)$.
Therefore
\begin{equation}\label{eq:3.3.3.3}
|u(x)|\leqslant c_p\phi(x)\text{ a.e. in }x\in\Omega.
\end{equation}
To prove that $|u(x)|\leqslant C_{p,\lambda}$ a.e. in $x\in\Omega$,
we let
\begin{equation}\label{eq:3.3.3.4}
w_2(x)=(u(x)-C_{p,\lambda})^+\ \text{ and }
\Omega_2=\{x\in\Omega:w_2(x) > 0\},
\end{equation}
then $w_2\in H^1_0(\Omega)$ and $u(x)|_{\Omega_2}\geqslant
C_{p,\lambda}>0$. Taking $v(x)=w_2(x)$ in (\ref{eq:3.1}) and using
(\ref{eq:3.3.3.3}) and $V(x)\geqslant1$, it follows that
\begin{equation}\nonumber
\int_{\Omega_2} \nabla u \nabla w_2+u  w_2dx\leqslant
 \int_{\Omega_2}|u|^{p-1}u  w_2-\frac{\lambda}{c_p}u^2  w_2dx\leqslant\int_{\Omega_2}
 C_{p,\lambda}  w_2dx,
\end{equation}
here we used the fact that
$\max\limits_{t\geqslant0}\{t^p-\frac{\lambda}{c_p}t^2\} =
C_{p,\lambda}$ if $p\in(1,2)$. This yields
\begin{equation}\nonumber
\int_{\Omega_2}|\nabla (u-C_{p,\lambda})^+ |^2 +
 |(u-C_{p,\lambda})^+|^2dx =\int_{\Omega_2}\nabla [u-C_{p,\lambda}] \nabla w_2 +
 [u-C_{p,\lambda}] w_2dx\leqslant0,
\end{equation}
it follows that $|\Omega_2|=0$, then $u(x)\leqslant C_{p,\lambda}$
a.e in $x\in\Omega$. Similarly,  $-u(x)\leqslant C_{p,\lambda}$.
Therefore
\begin{equation}\label{eq:3.3.3.5}
|u(x)|\leqslant C_{p,\lambda},\text{ a.e. in }x\in\Omega.
\end{equation}

{\bf Case 2:} $\Omega=\mathbb{R}^3$. In this case, it is not sure
if $w_1(x)$ given by (\ref{eq:3.2.1}) is in $H^1(\mathbb{R}^3)$,
so we need to replace $w_1(x)$ in (\ref{eq:3.2.1}) and the
correspondence by
\begin{equation}\nonumber
w_\epsilon(x)=(u(x)-c_{p} \phi(x)-\epsilon)^+ \ \text{ for any }
\epsilon>0.
\end{equation}
We claim that $w_\epsilon(x)\in H^1(\mathbb{R}^3)$ for any
$\epsilon>0$. In fact, by the second equation of (\ref{eq:3.1.L})
and Theorem 8.17 in \cite{DGilbargNStrudinger}, we know that
$$|\phi|_{L^\infty(B_1(y))}\leqslant C\left(|\phi|_{L^6(B_2(y))}+|u|^2_{L^{6}(B_3(y))}\right),\ \text{ for each }y\in\mathbb{R}^3,$$
where $C$ is a constant independent of $y\in\mathbb{R}^3$. This and
$\phi\in \mathscr{D}^{1,2}(\mathbb{R}^3)$ imply that $\phi
(x)\overset{|x|\rightarrow+\infty}{\rightarrow}0$.
 For any
$y\in\mathbb{R}^3$, taking $\varphi\in C^\infty_0(B_3(y))$, it
follows from (\ref{eq:3.1.L}) that
\begin{equation}\nonumber
\int_{B_3(y)}\nabla u \nabla\varphi+ b(x) u\varphi dx
=\int_{B_3(y)}c(x)\varphi dx,
\end{equation}
 where $b(x)=V(x)+\lambda\phi (x)$ and
$c(x)=|u|^{p-1}u(x)$. Clearly, $c(x)\in
L^{\frac{6}{p}}(\mathbb{R}^3)$ and $b(x)\in L^\infty(B^c_{R_1})$
for $R_1>0$ large enough since $\phi
(x)\overset{|x|\rightarrow+\infty}{\rightarrow}0$. For each
$|y|>R_1+3$, since $p<2$ and $\frac{6}{p} >\frac{3}{2}$, Theorem
8.17 in \cite{DGilbargNStrudinger} implies that
$$|u|_{L^\infty(B_1(y))}\leqslant C\left(|u|_{L^2(B_2(y))}+|c(x)|_{L^{\frac{6}{p}}(B_3(y))}\right),$$
where $C>0$ depends only on $p$ and
$|b(x)|_{L^\infty(B^c_{R_1})}$.
 This implies that
\begin{equation}\label{eq:3.3.1}
u(x)\rightarrow0 \text{ as } |x|\rightarrow\infty.
\end{equation}
 Hence,
there is $R_\epsilon>0$ such that  ${\rm supp}w_\epsilon\subset
B_{R_\epsilon}$, $w_\epsilon(x)\in H_0^1(B_{R_\epsilon})$ and
$w_\epsilon(x)\in H^1(\mathbb{R}^3)$. Noting that
$\nabla(u-c_p\phi)=\nabla(u-c_p\phi-\epsilon)$ for any
$\epsilon>0$, we see that (\ref{eq:3.3.3.2}) holds for
$w_\epsilon$, then Poincare inequality implies that
$\int_{\Omega_1} | w_\epsilon|^2dx \leq C_\epsilon \int_{\Omega_1}
|\nabla w_\epsilon|^2dx=
 0$ since $w_\epsilon(x)\in H_0^1(B_{R_\epsilon})$ and $w_\epsilon(x)\equiv
 0$ on $B_{R_\epsilon} \setminus \Omega_1$.
   So, $u(x)\leqslant c_p\phi(x)+\epsilon$ a.e. in $\mathbb{R}^3$
and (\ref{eq:3.3.3.3}) is obtained by letting
$\epsilon\rightarrow0$.
 By (\ref{eq:3.3.1}), there is $R_0>0$ such that ${\rm
supp}(u(x)-C_{p,\lambda})^+\subset B_{R_0}$ and $w_2(x)$ given by
(\ref{eq:3.3.3.4}) is in $H^1(\mathbb{R}^3)$, then exactly the
same
 as in Case 1 we get (\ref{eq:3.3.3.5}).  $\Box$

 To get the $L^\infty$-norm estimate of solutions to problem
(\ref{eq:2.1}) for $p\in[3,5)$, we need the following general
$L^\infty$-norm estimate for functions in
$\mathscr{D}^{1,2}(\mathbb{R}^N)$ by using its $L^{2^\ast}$ norm,
where $2^\ast=\frac{2N}{N-2}$.

\begin{lemma}\label{Le:3.2}
Let $N\geqslant3$, $p\in(1,\frac{N+2}{N-2})$ and let $u\in
\mathscr{D}^{1,2}(\mathbb{R}^3)\setminus\{0\}$ be a nonnegative
function such that
\begin{equation}\label{eq:3.5}
\int_{\mathbb{R}^{N}}\nabla u\nabla (h(u)\varphi)
dx\leqslant\int_{\mathbb{R}^{N}}|u|^{p-1}uh(u)\varphi dx,
\end{equation}
holds for any nonnegative $\varphi\in C_0^\infty(\mathbb{R}^N)$
and any nonnegative piecewise smooth function $h$ on $[0,+\infty)$
with $h(0)=0$ and $h'\in L^\infty(\mathbb{R}^3)$. Then, $u\in
L^\infty(\mathbb{R}^N)$ and there exist $C_1>0$ and $C_2>0$, which
depend only on $N$ and $p$, such that
\begin{eqnarray}\nonumber
|u|_{\infty}\leqslant
C_1\left(1+|u|_{2^\ast}^{C_2}\right)|u|_{2^\ast}\nonumber.
\end{eqnarray}
\end{lemma}
{\it Proof:} The main idea of the proof is Moser iterations, which
is somehow standard. For the sake of completeness, we give a proof
in the appendix based on \cite{TNS-ASNSPisa-1986} and
\cite{EgnellH-JDE-1992}. $\Box$
\begin{lemma}\label{Le:3.3} For each
$k\in\mathbb{N}$ large, let $p\in[3,5)$ and $u_k$ be a solution of
problem (\ref{eq:2.1}) given by Theorem \ref{th:2.1}. Then there
exist constants $M_0(\lambda)>0$ and $M(\lambda)>0$, independent
of $k$ and $\mu$, such that
\begin{equation}\nonumber
|u_k(x)|\leqslant M_0(\lambda)\ \text{a.e. in } x\in B_k\ \
\text{and }\  \|u_k\|_{k}+|\nabla \phi_k|_{L^{2}(B_k)}\leqslant
M(\lambda).
\end{equation}
Moreover, $M_0(\lambda)$ and $M(\lambda)$ are non-decreasing in
$\lambda>0$.
\end{lemma}
{\it Proof.} By Theorem \ref{th:2.1}, for any $k\in\mathbb{N}$
large, saying $k > k_0 \in \mathbb{N}$,  there exist $u_k$ and
$c_\lambda>0$ such that
\begin{equation}\label{eq:3.13.0}
I_{k}(u_k)=c_{\lambda,k} \leqslant c_\lambda\ \text{and }
I'_{k}(u_k)=0 \ \text{in } \mathscr{D}_k^\ast.
\end{equation}
 Similar to the
derivation of (\ref{eq:2.14.0}), we have
\begin{equation}\label{eq:3.13}
|\nabla u_k|_{L^2(B_k)}\leqslant \|u_k\|_{k}
\leqslant2\sqrt{c_{\lambda,k}},
\end{equation}
this and (\ref{eq:2.4}) yield
\begin{equation}\label{eq:3.14}
|\nabla \phi_k|_{L^{2}(B_k)}^2\leqslant C|
u_k|^4_{L^{\frac{12}{5}}(B_k)}\leqslant C\|u_k\|^4_{k}
\leqslant16Cc^2_{\lambda,k},
\end{equation}
where $C>0$ is independent of $k$, $\lambda$ and $\mu$. Let
$M(\lambda)=\sup\limits_{k>k_0,\mu>0}\{2\sqrt{c_{\lambda,k}}+4\sqrt{C}c_{\lambda,k}\}$,
then $M(\lambda)$ is finite for each $\lambda>0$ since
$c_{k,\lambda}\leqslant c_\lambda$. Moreover, $M(\lambda)$ is
non-decreasing in $\lambda$ since $c_{\lambda,k}$ is
non-decreasing in $\lambda$ by Theorem \ref{th:2.1}. Hence
(\ref{eq:3.13}) and (\ref{eq:3.14}) imply that $
\|u_k\|_{k}+|\nabla \phi_k|_{L^{2}(B_k)}\leqslant M(\lambda)$.\\
Define $\tilde{u}_k$ as in (\ref{eq:1.8}), then $\tilde{u}_k\in
\mathscr{D}_V$ and ${\rm supp}\tilde{u}_k\subseteq\bar{B_k}$,
hence (\ref{eq:3.13}) yields
\begin{equation}\label{eq:3.15}
 \|\tilde{u}_k\|_{\mathscr{D}_V}
\leqslant2\sqrt{c_{\lambda,k}}.
\end{equation}
For any $\varphi\in C^\infty_0(\mathbb{R}^3)$ and any nonnegative
piecewise smooth function $h$ with $h'\in L^\infty(\mathbb{R}^3)$
and $h(0)=0$, let $v=h(\tilde{u}_k^+)\varphi$, then $v \in
\mathscr{D}_k$ and $I_k'(u_k)v=0$, this yields
\begin{equation}\nonumber
\int_{\mathbb{R}^3}\nabla \tilde{u}^+_k\nabla v+V_\mu(x)
\tilde{u}^+_kvdx+\lambda\int_{\mathbb{R}^3}\tilde{\phi}_k
\tilde{u}^+_kvdx=\int_{\mathbb{R}^3}
|\tilde{u}^+_k|^{p-1}\tilde{u}^+_kvdx.
\end{equation}
Since $\tilde{\phi}_k$, $V_\mu$ and $\tilde{u}^+_k$ are nonnegative,
this implies that (\ref{eq:3.5}) holds for $u=\tilde{u}^+_k$ and
$N=3$. Then, by Lemma \ref{Le:3.2} and (\ref{eq:3.15}) there exist
$\bar{C}_1>0$ and $\bar{C}_2>0$, depend only on $p$, such that
\begin{equation}\nonumber
|\tilde{u}^+_k|_{{\infty}}\leqslant
C_1(1+|\tilde{u}_k^+|_{2^\ast}^{C_2})|\tilde{u}_k^+|_{2^\ast}
\leqslant\sup\limits_{k>k_0,\mu>0}\{\bar{C}_1(1+c^{\bar{C}_2}_{\lambda,k})
\sqrt c_{\lambda,k}\}:=M_0(\lambda).
\end{equation}
 Similarly, we know that $|\tilde{u}^-_k|_{{\infty}}\leqslant
M_0(\lambda)$. So, $|{u}_k|_{L^{\infty}(B_k)} \leq
|\tilde{u}_k|_{{\infty}}\leqslant M_0(\lambda) $. Moreover,
$M_0(\lambda)$ is nondecreasing in $\lambda>0$ since
$c_{\lambda,k}$ is nondecreasing in $\lambda$.    $\Box$\\

\section{Proofs of Theorems \ref{th:1.1} and \ref{th:1.2}}
For $r>0$, let $\xi_r\in C^\infty(\mathbb{R}^3)$ such that
\begin{equation}\label{eq:4.4}
\xi_r(x)=\left\{\begin{array}{ll}
 1,\ \ \  |x|>\frac{r}{2}, \\
 0,\ \ \ |x|<\frac{r}{4},
 \end{array}\right.\ \ \text{ with }\  |\nabla\xi_r|\leqslant \frac{8}{r}.
\end{equation}
To prove our Theorems, we need the following two lemmas.
\begin{lemma}\label{Le:4.1} Assume $\rm(G_1)$ $\rm (G_3)$ hold. Let $p\in(1,2)$, $u\in
\mathscr{D}_V$ and $\phi\in L^6(\mathbb{R}^3)$ be such that
$|u|_\infty\leqslant L$ and $ |u(x)|\leqslant K\phi(x)\ \text{a.e.
in }x\in\mathbb{R}^3$, for some $L>0$ and $K>0$. Moreover, $\text{
for all }\eta\in C^1_B(\mathbb{R}^3)\text{ with }\eta\geqslant0$,
there holds
\begin{equation}\label{eq:4.0.0}
\int_{\mathbb{R}^3}\nabla u\nabla (u\eta)+V_\mu(x) u^2\eta
dx+\lambda\int_{\mathbb{R}^3}\phi u^2\eta dx=\int_{\mathbb{R}^3}
|u|^{p+1}\eta dx,
\end{equation}
where $\lambda>0$. Then,
 if $\mu>\bar{\mu}=\max\{0,L^{p-1}-1\}$, there exists
 $M_\mu>0$ (depends on $L,K,p,\lambda,\mu$) which is decreasing in $\mu>\bar{\mu}$ such that
$$\|  u\|_{\mathscr{D}_V}+\left\{\int_{\mathbb{R}^3}\phi u^2dx\right\}^{\frac{1}{2}}\leqslant M_\mu.$$
\end{lemma}
{\it Proof:} For $\xi_r$ given by (\ref{eq:4.4}), taking
 $\eta=\xi_r$ in (\ref{eq:4.0.0}), it gives that
\begin{equation}\nonumber
\int_{\mathbb{R}^3}(|\nabla
u|^2+V_\mu(x)u^2)\xi_rdx+\lambda\int_{\mathbb{R}^3}\phi
u^2\xi_rdx=\int_{\mathbb{R}^3} |u|^{p+1}\xi_rdx-\int_{\mathbb{R}^3}
u\nabla u\nabla\xi_rdx,
\end{equation}
then, by $|u|_\infty\leqslant L$ and (\ref{eq:4.4}) we have
\begin{equation}\nonumber
\begin{split}
&\int_{\mathbb{R}^3}(|\nabla u|^2+(1+\mu
g(x))u^2)\xi_rdx+\lambda\int_{\mathbb{R}^3}\phi
u^2\xi_rdx\\
&\leqslant L^{p-1}\int_{\mathbb{R}^3}
 u^{2}\xi_rdx+\frac{8}{r}\int_{\mathbb{R}^3}(|\nabla
u|^2+u^2)dx,
\end{split}
\end{equation}
that is,
\begin{equation}\label{eq:4.5}
\int_{\mathbb{R}^3}(|\nabla u|^2+[1+\mu g(x)-L^{p-1}]u^2)\xi_rdx
\leqslant \frac{8}{r}\int_{\mathbb{R}^3}(|\nabla u|^2+u^2)dx.
\end{equation}
Since $\mu>\bar{\mu} \geq L^{p-1}-1$, for each
$\tau\in(\bar{\mu},\mu)$, there exists $\epsilon_\tau>0$ such that
$\tau(1-\epsilon_\tau)>L^{p-1}-1$. By $\rm(G_3)$ we can find
$R_\tau>0$ such that $ g(x)\geqslant 1-\epsilon_\tau$ for all
$|x|\geqslant R_\tau$, then $1+\mu
g(x)>1+\tau(1-\epsilon_\tau)>L^{p-1}$ for all $|x|\geq R_\tau$ and
there is $\delta_\tau>0$ such that $1+\mu g(x)>
L^{p-1}+\delta_\tau$ for all $|x|\geqslant R_\tau$. So, if
$\mu>\bar{\mu}$ and $r/4>R_\tau$, it follows from (\ref{eq:4.5})
that, for each $\tau \in (\bar{\mu},\mu)$,
\begin{equation}\label{eq:4.6}
\int_{|x|>\frac{r}{2}}(|\nabla u|^2+\delta_\tau u^2)dx \leqslant
\frac{8}{r}\int_{\mathbb{R}^3}(|\nabla u|^2+u^2)dx.
\end{equation}
On the other hand,   taking $\eta\equiv1$ in (\ref{eq:4.0.0}), it
gives
\begin{equation}\label{eq:4.7}
\int_{\mathbb{R}^3}|\nabla u|^2+(1+\mu
g(x))u^2dx\leqslant\int_{\mathbb{R}^3}L^{p-1}
 u^{2}-\frac{\lambda}{K}  |u|^{3}dx.
\end{equation}
Let $\beta(t)=L^{p-1}t^2-\frac{\lambda}{K}t^3$($t\geqslant0$), then
there exists $C^\ast:=C(L,K,p,\lambda)\in(0,+\infty)$ such that
$\beta(t)\leqslant C^\ast<+\infty$ for all $t\geqslant0$. Therefore,
if $r>4R_\tau$ and $\mu>\bar{\mu}$, it follows from (\ref{eq:4.6})
and (\ref{eq:4.7}) that
\begin{equation}\nonumber
\begin{split}
&\int_{\mathbb{R}^3}|\nabla u|^2+(1+\mu g(x))u^2dx
\leqslant\int_{|x|<\frac{r}{2}}\beta(u)dx+\int_{|x|\geqslant\frac{r}{2}}\beta(u)dx\\
&\leqslant\int_{|x|<\frac{r}{2}}C^\ast
dx+L^{p-1}\int_{|x|\geqslant\frac{r}{2}}u^2dx \leqslant
C^\ast|B_{\frac{r}{2}}(0)|
+L^{p-1}\frac{8}{r\delta_\tau}\int_{\mathbb{R}^3}(|\nabla
u|^2+u^2)dx.
\end{split}
\end{equation}
Take $r_\tau>4R_\tau$ large enough such that
$L^{p-1}\frac{8}{r_\tau\delta_\tau}<\frac{1}{2}$, then
\begin{equation}\nonumber
\begin{split}
\int_{\mathbb{R}^3}|\nabla u|^2+&(1+\mu
g(x))u^2dx\leqslant2C^\ast|B_{\frac{r_\tau}{2}}|.
\end{split}
\end{equation}
By Sobolev embedding, $|u|_{p+1}$ is also bounded above by a
constant depending  on $C^*$ and $\tau$. This and (\ref{eq:4.0.0})
imply that there exists $M_\tau:=M(L,K,p,\lambda,\tau)$ such that
$$\|u\|_{\mathscr{D}_V}+\left\{\int_{\mathbb{R}^3}\phi u^2dx\right\}^{\frac{1}{2}}\leqslant M_\tau,
 \mbox{ for each } \tau \in (\bar{\mu},\mu).$$
Let $M_\mu=\inf\limits_{\tau\in(\bar{\mu},\mu)}M_\tau$, then
$M_\mu$ is decreasing in $\mu>\bar{\mu}$ and depends only on
$L,K,p,\lambda,\mu$ such that
$$\|u\|_{\mathscr{D}_V}+\left\{\int_{\mathbb{R}^3}\phi u^2dx\right\}^{\frac{1}{2}}\leqslant M_\mu. \ \ \ \ \ \ \ \ \ \  \Box$$
\begin{lemma}\label{Le:4.2} If $p>1$, and
$\rm(G_1)$ $\rm(G_3)$ hold, let $\{ u_k\}$ be bounded in $
\mathscr{D}_V$ such that $|u_k|_\infty\leqslant \bar{M}$ (for some
$\bar{M}>0$) and
\begin{equation}\label{eq:4.8}
\int_{\mathbb{R}^3}\nabla u_k\nabla (u_k\eta)+V_\mu(x) u^2_k\eta
dx\leqslant\int_{\mathbb{R}^3}  |u_k|^{p+1}\eta dx,
\end{equation}
 for all $\eta\in
C_B^1(\mathbb{ R}^3) \text{ with } \eta\geqslant0$. Then, for each
$\mu>\max\{0,\bar{M}^{p-1}-1\}$, there exists $u\in \mathscr{D}_V$
such that,  passing to a subsequence,
\begin{equation}\label{eq:4.8.0}
  |  u_k-u|_q\overset{k}{\rightarrow}0\
\text{ for }\  q\in[2,6).
\end{equation}
\end{lemma}
{\it Proof:} Since $\{u_k\}$ is bounded in $\mathscr{D}_V$, there
exists $u\in \mathscr{D}_V$ such that, passing to a subsequence,
\begin{equation}\label{eq:4.10}
u_k\overset{k}{\rightharpoonup}u \ \text{ weakly in }
\mathscr{D}_V,\ \ u_k(x)\overset{k}{\rightarrow}u(x)\ \text{a.e. in
} x\in\mathbb{R}^3.
\end{equation}
 For $\xi_r$ given by (\ref{eq:4.4}), taking
$\eta=\xi_r$ in (\ref{eq:4.8}),
\begin{equation}\nonumber
\int_{\mathbb{R}^3}\nabla u_k\nabla (
 u_k\xi_r)+V_\mu(x)u_k^2\xi_rdx\leqslant\int_{\mathbb{R}^3}
 |u_k|^{p+1}\xi_rdx.
\end{equation}
 Since
$\mu>\max\{0,\bar{M}^{p-1}-1\}$, $|u_k|_\infty\leqslant\bar{ M}$
and $\rm(G_3)$ holds, similar to the discussion of (\ref{eq:4.6}),
there exists $R_\mu>0$ and $\delta_\mu>0$ such that, for all
$R\geqslant4R_\mu$, we have
\begin{equation}\label{eq:4.11}
\int_{|x|>R}(|\nabla u_k|^2+\delta_\mu u_k^2)dx <\frac{T}{R}\ \
 \text{ uniformly for } k\in\mathbb{N}\ \text{large},
\end{equation}
where $T=\sup\limits_{k\in\mathbb{N}}\|u_k\|_{\mathscr{D}_V}$. Since
$H^1(B_{R})\hookrightarrow L^q(B_{R})$ is compact for $1\leqslant
q<6$, passing to a subsequence, (\ref{eq:4.10}) implies that
\begin{equation}\label{eq:4.12}
u_k(x)\overset{k}{\rightarrow}u(x)\ \text{ in }L^q(B_{R}) \text{ for
}  1\leqslant q<6.
\end{equation}
For $k\in\mathbb{N}$ large and any $R>4R_\mu$ large enough, we
have
\begin{equation}\nonumber
\begin{split}
 &|  u_k-u|_q\leqslant|  u_k-u|_{L^q(B_{R})}+|  u_k-u|_{L^q(B^c_{R})}\\
 &\leqslant|  u_k-u|_{L^q(B_{R})}
 +|  u_k-u|^{\frac{6-q}{2}}_{L^2(B^c_{R})}|  u_k-u|^{\frac{3q-6}{2}}_{L^6(B^c_{R})}\\
 &\leqslant|  u_k-u|_{L^q(B_{R})}
 +C\|  u_k-u\|^{\frac{3q-6}{2}}_{H^1(\mathbb{R}^3)}\left(|u|^{\frac{6-q}{2}}_{L^2(B^c_{R})}
 +|  u_k|^{\frac{6-q}{2}}_{L^2(B^c_{R})}\right)\\
 &\leqslant|  u_k-u|_{L^q(B_{R})}
 +C\left(|u|^{\frac{6-q}{2}}_{L^2(B^c_{R})}
 +{(T/R)}^{\frac{6-q}{2}}\right)\text{ by
(\ref{eq:4.11})}.
\end{split}
\end{equation}
By letting $k\rightarrow+\infty$, then $R\rightarrow+\infty$, we get
(\ref{eq:4.8.0}). \ \ \ $\Box $

\begin{lemma}\label{Le:4.3} Assume $p\in(1,5)$.
If $\rm(G_1)$ $\rm(G_3)$ hold and $\{ u_k\}$ is bounded in $
\mathscr{D}_V$ satisfying
\begin{equation}\nonumber
|u_k|_\infty\leqslant\bar{M},\ \text{ for some } \bar{M}>0.
\end{equation}
 Let $\phi_k$
be the solution of $-\Delta\phi=u_k^2$ in $\mathscr{D}_0^{1,2}(B_k)$
and $\tilde{\phi}_k$ be the extension of $\phi_k$ on $\mathbb{R}^3$
defined as (\ref{eq:1.8}) such that, for all $\eta\in C_B^1(\mathbb{
R}^3) \text{ with } \eta\geqslant0$,
\begin{equation}\label{eq:4.30.1}
\int_{\mathbb{R}^3}\nabla u_k\nabla (u_k\eta)+V_\mu(x) u^2_k\eta
dx+\lambda\int_{\mathbb{R}^3}\tilde{\phi}_ku_k^2\eta
dx=\int_{\mathbb{R}^3} |u_k|^{p+1}\eta dx,
\end{equation}
where $\lambda>0$. Then, for each $\mu>\max\{0,\bar{M}^{p-1}-1\}$,
there exists $u\in \mathscr{D}_V$ such that,  passing to a
subsequence,
\begin{equation}\label{eq:4.30.3}
\|  u_k-u\|_{\mathscr{D}_V}\overset{k}{\rightarrow}0 \text{ and }\
|\nabla(\tilde{\phi}_k-\phi)|_2\overset{k}{\rightarrow}0\text{ with
} \phi=\frac{1}{4\pi}\int_{\mathbb{R}^3}\frac{u^2(y)}{|x-y|}dy.
\end{equation}
\end{lemma}
{\it Proof: } Since $\lambda>0$ and $\tilde{\phi}_k\geqslant0$,
(\ref{eq:4.30.1}) implies (\ref{eq:4.8}), for each
$\mu>\max\{0,\bar{M}^{p-1}-1\}$  applying Lemma \ref{Le:4.2} we
have $u\in \mathscr{D}_V$ such that,
\begin{equation}\label{eq:4.30.4}
  |  u_k-u|_q\overset{k}{\rightarrow}0\
\text{ for }\  q\in[2,6),
\end{equation}
This and Lemma \ref{Le:2.3} imply that, passing to a subsequence,
\begin{equation}\nonumber
  \tilde{\phi}_k\overset{k}{\rightarrow}
   \phi=\frac{1}{4\pi}\int_{\mathbb{R}^3}\frac{u^2(y)}{|x-y|}dy \ \text{ strongly in } \mathscr{D}_0^{1,2}(\mathbb{R}^3).
\end{equation}
Hence,
\begin{equation}\label{eq:4.30.5}
 \int_{\mathbb{R}^3} \tilde{\phi}_k u_k^2dx\overset{k}{\rightarrow}
   \int_{\mathbb{R}^3}\phi u^2dx.
\end{equation}
Then (\ref{eq:4.30.3}) follows from (\ref{eq:4.30.4}),
(\ref{eq:4.30.5}) and (\ref{eq:4.30.1}) with $\eta\equiv1$.\ \ \
$\Box$ \\

 Now, we are ready to prove Theorems \ref{th:1.1} and
\ref{th:1.2}. For $u_k\in \mathscr{D}_k$, $\phi_k\in
\mathscr{D}^{1,2}_0(B_k)$ and $\alpha>0$, $c_\lambda>0$ given by
Theorem \ref{th:2.1}, let $\tilde{u}_k\in \mathscr{D}_V$ and
$\tilde{\phi}_k\in \mathscr{D}^{1,2}(\mathbb{R}^3)$ be the
extensions of $u_k$ and $\phi_k$ on $\mathbb{R}^3$ defined as
(\ref{eq:1.8}) respectively.
Then,
\begin{equation}\label{eq:4.1}
\frac{1}{2}\int_{\mathbb{R}^3}|\nabla
\tilde{u}_k|^2+V_\mu(x)\tilde{u}_k^2
dx+\frac{\lambda}{4}\int_{\mathbb{R}^3}\tilde{\phi}_{k}(x)\tilde{u}_k^2dx-\frac{1}{p+1}\int_{\mathbb{R}^3}
 |\tilde{u}_k|^{p+1}dx=c_{\lambda,k},
\end{equation}
and $c_{\lambda,k}\in[\alpha,c_\lambda]$. Moreover, for any
$\varphi\in \mathscr{D}_V$ with ${\rm supp}\varphi\subseteq
\bar{B}_k$, by $I'_k(u_k)=0$ we see that
\begin{equation}\label{eq:4.2}
\int_{\mathbb{R}^3}\nabla \tilde{u}_k\nabla\varphi+V_\mu(x)
\tilde{u}_k\varphi dx+\lambda\int_{\mathbb{R}^3}\tilde{\phi}_{k}
\tilde{u}_k\varphi dx=\int_{\mathbb{R}^3}
|\tilde{u}_k|^{p-1}\tilde{u}_k\varphi dx.
\end{equation}
Note that for any $\varphi \in C^\infty_0(\mathbb{R}^3)$,
$\text{supp}\varphi \subset \bar{B}_k$ as $k\rightarrow+\infty$,
therefore,
\begin{equation}\label{eq:4.3}
\int_{\mathbb{R}^3}\nabla \tilde{u}_k\nabla\varphi+V_\mu(x)
\tilde{u}_k\varphi dx+\lambda\int_{\mathbb{R}^3}\tilde{\phi}_{k}
\tilde{u}_k\varphi dx=\int_{\mathbb{R}^3}| \tilde{u}_k|^{p-1}
\tilde{u}_k\varphi dx+o(1).
\end{equation}
{\bf Proof of Theorem \ref{th:1.1}:} Let $\lambda_\ast$ be given by
Theorem \ref{th:2.1}, and $\lambda_\ast<+\infty$ if
 $p\in(1,2)$. For $\lambda\in(0,\lambda_\ast)$, let $
(u_k,\phi_k)\in\mathscr{D}_k\times\mathscr{D}_0^{1,2}(B_k)$ be a
solution of (\ref{eq:2.1}) for each $k\in\mathbb{N}$ large. Applying
Lemma \ref{Le:3.1} with $\Omega=B_k$, there exist $c_p>0$ and
$C_{p,\lambda}>0$ such that
\begin{equation}\nonumber
|u_k(x)|\leqslant c_p\phi_k(x)\ \text{and } \ |u_k(x)|\leqslant
C_{p,\lambda},\ \
  \text{a.e. in } x\in B_k.
\end{equation}
Then,
\begin{equation}\nonumber
|\tilde{u}_k(x)|\leqslant c_p\tilde{\phi}_k(x)\ \text{and } \
|\tilde{u}_k(x)|\leqslant C_{p,\lambda},\ \
  \text{a.e. in } x\in \mathbb{R}^3.
\end{equation}
Since $\tilde{u}_k\eta\in\mathscr{D}_k$ for any $\eta\in
C^1_B(\mathbb{R}^3)$, it follows from $I'_k(u_k)=0$ that
(\ref{eq:4.0.0}) holds for $u=\tilde{u}_k$ and
$\phi=\tilde{\phi}_k$, applying Lemma \ref{Le:4.1} with
$L=C_{p,\lambda}$, $K=c_p$ and $u=\tilde{u}_k$, there is
$M_{\mu}:=M(p,\lambda,\mu)>0$, which is decreasing in $\mu>\mu_1$,
such that
$$\|  \tilde{u}_k\|_{\mathscr{D}_V}+|\nabla\tilde{\phi}_k|_2=\|  \tilde{u}_k\|_{\mathscr{D}_V}+\left\{\int_{\mathbb{R}^3}{\tilde{\phi}}_k
\tilde{u}_k^2dx\right\}^{\frac{1}{2}}\leqslant M_{\mu}. $$
Similarly, $\{\tilde{u}_k\}$ and $\{\tilde{\phi}_k\}$ satisfy
(\ref{eq:4.30.1}) by $I'_k(u_k)=0$. Hence, Lemma \ref{Le:4.3} shows
that there exists $u\in
 \mathscr{D}_V$
such that, passing to a subsequence
\begin{equation}\label{eq:4.16}
\tilde{u}_k\overset{k}{\rightarrow}u\ \text{ strongly in }
\mathscr{D}_V \ \text{ and }\
\tilde{\phi}_k\overset{k}{\rightarrow}\phi_u\ \text{ strongly in }
\mathscr{D}^{1,2}(\mathbb{R}^3).
\end{equation}
Furthermore,
\begin{equation}\label{eq:4.15}
\|u\|_{\mathscr{D}_V}+|\nabla\phi_{u} |_{2}\leqslant M_{\mu}\ \text{
and } \ |u|_\infty\leqslant C_{p,\lambda}.
\end{equation}
Combining (\ref{eq:4.1}) and (\ref{eq:4.16}), we have $
I(u)\in[\alpha,c_\lambda] $. On the other hand, by (\ref{eq:4.3}),
for any $\varphi\in C^\infty_0(\mathbb{R}^3)$, as
$k\rightarrow+\infty$, $ I'(\tilde{u}_k)\varphi=o(1)$, then, using
(\ref{eq:4.16}) we see that $ I'(u)\varphi=0$. So, $u$ is a
nontrivial solution of (\ref{eq:1.1}) in $ \mathscr{D}_V$. By the
second equation of (\ref{eq:1.1}) and Theorem 8.17 in
\cite{DGilbargNStrudinger}, for any $y\in \mathbb{R}^3$, we have
\begin{equation}\nonumber
\sup\limits_{B_1(y)}\phi_{u}\leqslant
C(|\phi_{u}|_{L^6(B_2(y))}+|u|^2_{L^6(B_4(y))}),
\end{equation}
where $C>0$ is a constant independent of $u$ and $y\in
\mathbb{R}^3$. This and (\ref{eq:4.15}) show that
\begin{equation}\label{eq:4.17.0}
\sup\limits_{x\in\mathbb{R}^3}\phi_{u}(x)\leqslant
CM_{\mu}(1+M_{\mu})\triangleq C(p,\lambda,\mu)=:C_\mu,
\end{equation}
where $C_\mu>0$ is a constant dependents only on $p,\lambda,\mu$.
$\Box$\\
 {\bf Proof of Theorem \ref{th:1.2}:} For $p\in[3,5)$ and $\lambda>0$,
 let $(u_k,\phi_k)\in
\mathscr{D}_k\times\mathscr{D}^{1,2}(B_k)$ be a solution given by
Theorem \ref{th:2.1}. By
 Lemma \ref{Le:3.3}, we know that
 $$|u_k|_\infty\leqslant M_0(\lambda),\ \ \ \|  u_k\|_{k}+|\nabla\phi_k
|_{L^2(B_k)}\leqslant M(\lambda), $$
 where $M_0(\lambda)$ and
$M(\lambda)$ are given by Lemma \ref{Le:3.3}. Hence,
$$|\tilde{u}_k|_\infty\leqslant M_0(\lambda),\ \ \ \|  \tilde{u}_k\|_{\mathscr{D}_V}+|\nabla\tilde{\phi}_k
|_{2}\leqslant M(\lambda). $$ Applying Lemma \ref{Le:4.3} with
$u_k=\tilde{u}_k$, there exists $u\in
 \mathscr{D}_V$
such that,
\begin{equation}\label{eq:4.19}
\tilde{u}_k\overset{k}{\rightarrow}u \ \text{  in } \mathscr{D}_V.
\end{equation}
Moreover,
\begin{equation}\label{eq:4.18}
\|u\|_{\mathscr{D}_V}+|\nabla\phi_{u} |_{2}\leqslant M(\lambda)\ \
\text{ and } \ \ |u|_\infty\leqslant M_0(\lambda).
\end{equation}
Hence, (\ref{eq:4.19}) and (\ref{eq:4.1}) show that $
I(u)\in[\alpha,c_\lambda]$, where $\alpha$ and $c_\lambda$ are
constants given by Theorem \ref{th:2.1}.
 On the other hand, by
(\ref{eq:4.3}), for any $\varphi\in C^\infty_0(\mathbb{R}^3)$, as
$k\rightarrow+\infty$, $ I'(  \tilde{u}_k)\varphi=o(1) $, and it
follows from (\ref{eq:4.19}) that $ I'(u)\varphi=0 $. So, $u$ is a
nontrivial solution of (\ref{eq:1.1}) in $ \mathscr{D}_V$. Noting
(\ref{eq:4.18}), similar to the discussion of (\ref{eq:4.17.0}), we
have
\begin{equation}\nonumber
\sup\limits_{x\in\mathbb{R}^3}\phi_{u}(x)\leqslant
CM(\lambda)(1+M(\lambda))=:M_1(\lambda),
\end{equation}
where $M_1(\lambda)$ is a constant and independent of $\mu$. By
Lemma \ref{Le:3.3}, $M(\lambda)$ is non-decreasing in $\lambda>0$,
hence $M_1(\lambda)$ is
also non-decreasing in $\lambda$. $\Box$ \\

\section{Proofs of Theorems \ref{th:1.3} to \ref{th:1.5}}
{\bf Proof of Theorem \ref{th:1.3}:} By Theorem \ref{th:1.2}, if
$\mu>\mu_2(\lambda)$, problem (\ref{eq:1.1}) has always a solution
$u_\lambda\in\mathscr{D}_V$ for any $\lambda\in(0,1)$ such that
$$\|u_\lambda\|_{\mathscr{D}_V}+|\nabla\phi_{u_\lambda} |_{2}\leqslant M(\lambda)\leqslant M(1)\
\text{ and } \ \ |u_\lambda|_\infty\leqslant M_0(\lambda)\leqslant
M_0(1),$$ since $M_0(\lambda)$ and $M(\lambda)$ are nondecreasing in
$\lambda>0$. Then, there exists $u_0\in \mathscr{D}_V$, passing to a
subsequence, such that
\begin{equation}\label{eq:5.0}
u_\lambda\overset{\lambda\rightarrow0}{\rightharpoonup}u_0 \text{
weakly in } \mathscr{D}_V \text{ and }
u_\lambda\overset{\lambda\rightarrow0}{\rightarrow}u_0, \text{ a.e.
in }x\in\mathbb{R}^3.
\end{equation}
 For any $\varphi\in \mathscr{D}_V$,
$ \int_{\mathbb{R}^3}\phi_{u_\lambda} u_\lambda \varphi
dx\leqslant|\phi_{u_\lambda}|_6|u_\lambda|_{\frac{12}{5}}|\varphi|_{\frac{12}{5}}\leqslant
C(M(1))$ and
\begin{equation}\label{eq:4.21}
\int_{\mathbb{R}^3}\nabla u_\lambda\nabla \varphi+V_\mu(x)u_\lambda
\varphi dx+\lambda\int_{\mathbb{R}^3}\phi_{u_\lambda} u_\lambda
\varphi dx=\int_{\mathbb{R}^3}|u_\lambda|^{p-1}u_\lambda\varphi dx,
\end{equation}
let $\lambda\rightarrow0$, we have
\begin{equation}\label{eq:4.21.0}
\int_{\mathbb{R}^3}\nabla u_0\nabla \varphi+V_\mu(x)u_0 \varphi
dx=\int_{\mathbb{R}^3}|u_0|^{p-1}u_0 \varphi dx,
\end{equation}
that is, $u_0$ is a weak solution of (\ref{eq:1.1}) with
$\lambda=0$. Since $u_\lambda \eta\in\mathscr{D}_V$ for any
$\eta\in C_B^1(\mathbb{R}^3)\text{ with }\eta\geqslant0$, by
(\ref{eq:4.21}) and $\lambda>0$, we see that
\begin{equation}\nonumber
\int_{\mathbb{R}^3}\nabla u_\lambda\nabla
(u_\lambda\eta)+V_\mu(x)u^2_\lambda \eta
dx\leqslant\int_{\mathbb{R}^3}|u_\lambda|^{p+1} \eta dx.
\end{equation}
By Lemma \ref{Le:4.2}, there is $u^\ast\in\mathscr{D}_V$ such that
$u_\lambda\overset{\lambda\rightarrow0}{\rightarrow}u^\ast$ in
$L^q(\mathbb{R}^3)$ for $q\in[2,6)$. Hence $u^\ast=u_0$ by
(\ref{eq:5.0}). Let $\varphi=u_\lambda$ in (\ref{eq:4.21}) and
$\varphi=u_0$ in (\ref{eq:4.21.0}), since
$\int_{\mathbb{R}^3}\phi_{u_\lambda} u^2_\lambda dx=|\nabla
\phi_{u_\lambda}|^2_2\leqslant M^2(1)$, as $\lambda\rightarrow0$, we
have $u_\lambda\rightarrow u_0$ in $\mathscr{D}_V$. This and
$I(u_\lambda)\in[\alpha,c_\lambda]$ for $\lambda>0$ imply that
$I_0(u_0)\geqslant\alpha>0$, hence $u_0\not\equiv0$. $\Box$\\
 {\bf Proof of Theorem
\ref{th:1.4}.} To make it clear that $\mathscr{D}_V$ depends on
$\mu$, here we denote $\mathscr{D}_V$ by $\mathscr{D}_\mu$. Let
$\mu_0>\max\{C_{p,\lambda}^{p-1}-1, M_0^{p-1}-1\}$ and $\mu_0>0$
by Remark \ref{r:1.1}. For $\mu\geqslant\mu_0$, let $u_\mu$ be the
solution given by Theorem \ref{th:1.1} or Theorem \ref{th:1.2},
then $\{u_\mu\}$ is bounded in $D_{\mu_0}$ since $M_\mu \leq
M_{\mu_0}$.
 For some $\bar{u}\in
\mathscr{D}_{\mu_0}\subset{H^1(\mathbb{R}^3)}$ we may assume that
\begin{equation}\label{eq:5.6}
u_\mu\rightharpoonup \bar{u}\ \ \text{weakly in }
\mathscr{D}_{\mu_0} \text{ and } u_\mu\rightarrow \bar{u} \ a.e. \
\text{in} \ x\in\mathbb{R}^3, \text{ as }\mu\rightarrow+\infty.
\end{equation}
Then, Fatou's lemma shows that
\begin{equation}\nonumber
\int_{\mathbb{R}^3} g(x)\bar{u}^2
dx\leqslant\lim\limits_{\mu\rightarrow+\infty}\int_{\mathbb{R}^3}
g(x)u_\mu^2dx \leqslant
\lim\limits_{\mu\rightarrow+\infty}\frac{1}{\mu}\|u_\mu\|^2_{\mathscr{D}_{\mu}}\leqslant
\lim\limits_{\mu\rightarrow+\infty}\frac{M^2_{\mu_0}}{\mu}=0,
\end{equation}
so $\rm(G_1)$ and $\rm(G_2)$ means that $\bar{u}(x)=0$ a.e. in
$x\in\mathbb{R}^3\setminus \Omega_0$. On the other hand, for any
$\varphi\in \mathscr{D}_{\mu}$, by $I'(u_\mu)=0$ we have
\begin{equation}\label{eq:5.5}
 \int_{\mathbb{R}^3}\nabla u_\mu\nabla\varphi +(1+\mu g(x))u_\mu\varphi
 dx+\lambda\int_{\mathbb{R}^3}\phi_{u_\mu}u_\mu\varphi dx
  =\int_{\mathbb{R}^3} |u_\mu|^{p-1}u_\mu\varphi dx.
\end{equation}
 For each $\mu\geqslant\mu_0$ and
$\eta\in C^1_B(\mathbb{R}^3)$ with $\eta\geqslant0$, it follows from
(\ref{eq:5.5}) and $\rm(G_1)$ that
\begin{equation}\nonumber
 \int_{\mathbb{R}^3}\nabla u_\mu\nabla (u_\mu\eta) +(1+\mu_0 g(x))u^2_\mu\eta dx\leqslant\int_{\mathbb{R}^3} |u_\mu|^{p+1}\eta
 dx.
\end{equation}
Then, it follows from (\ref{eq:5.6}), Theorem \ref{th:1.1} or
\ref{th:1.2} and Lemma \ref{Le:4.2} with
$\mathscr{D}_V=\mathscr{D}_{\mu_0}$ and $\{u_k\}=\{u_\mu\}$ that
\begin{equation}\label{eq:5.5.3}
  |u_\mu-\bar{u}|_q\overset{\mu\rightarrow+\infty}{\rightarrow}0 \text{ for }
q\in[2,6),
\end{equation}
hence, Lemma \ref{Le:2.2} yields that
\begin{equation}\label{eq:5.5.3.0}
 \int_{\mathbb{R}^3}
 \phi_{u_\mu}u^2_\mu(x)dx\overset
 {\mu\rightarrow+\infty}{\rightarrow}\int_{\mathbb{R}^3}
 \phi_{\bar{u}}\bar{u}^2(x)dx.
\end{equation}
Let $\varphi\in C^\infty_0(\mathbb{R}^3)$, by (\ref{eq:5.6}) and
(\ref{eq:5.5.3.0}) (see (3.18) in \cite{JZ-ActMS} for the
details), we have
\begin{equation}\label{eq:5.5.3.1}
\int_{\mathbb{R}^3}
 \phi_{u_\mu}u_\mu(x)\varphi(x)dx\overset
 {\mu\rightarrow+\infty}{\rightarrow}\int_{\mathbb{R}^3}
 \phi_{\bar{u}}\bar{u}(x)\varphi(x)dx,
\end{equation}
and  (\ref{eq:5.5.3.1})  holds also for $\varphi \in
H^1(\mathbb{R}^3)$ by the density of $C^\infty_0(\mathbb{R}^3)$ in
$H^1(\mathbb{R}^3)$. Let $\varphi=\bar{u}$ in (\ref{eq:5.5}) and
in (\ref{eq:5.5.3.1}), then, by (\ref{eq:5.6}) and
(\ref{eq:5.5.3}) we have
\begin{equation}\label{eq:5.5.2}
 \int_{\mathbb{R}^3}|\nabla \bar{u}|^2 +\bar{u}^2 dx+
 \lambda\int_{\mathbb{R}^3}\phi_{\bar{u}}\bar{u}^2dx=\int_{\mathbb{R}^3}|\bar{u}|^{p+1}
 dx.
\end{equation}
By (\ref{eq:5.5.3}) and (\ref{eq:5.5}) with $\varphi=u_\mu$, we see
that
\begin{equation}\label{eq:5.5.1}
 \int_{\mathbb{R}^3}|\nabla u_\mu|^2 +u^2_\mu dx+\lambda\int_{\mathbb{R}^3}\phi_{u_\mu}
u_\mu^2(x)dx\leqslant\int_{\mathbb{R}^3} |u_\mu|^{p+1}
 dx\overset{\mu\rightarrow+\infty}{\rightarrow}\int_{\mathbb{R}^3}|\bar{u}|^{p+1}
 dx.
\end{equation}
It follows from (\ref{eq:5.5.2}) and (\ref{eq:5.5.1}) that
\begin{equation}\nonumber
 \lim\limits_{\mu\rightarrow+\infty}\int_{\mathbb{R}^3}|\nabla
u_\mu|^2 +u^2_\mu dx
 +\lambda\int_{\mathbb{R}^3}\phi_{u_\mu}u^2_\mu dx
 \leqslant\int_{\mathbb{R}^3}|\nabla \bar{u}|^2 +\bar{u}^2
dx+\lambda\int_{\mathbb{R}^3}\phi_{\bar{u}}\bar{u}^2dx,
\end{equation}
this and (\ref{eq:5.5.3.0}) show that
\begin{equation}\nonumber
 \lim\limits_{\mu\rightarrow+\infty}\int_{\mathbb{R}^3}|\nabla u_\mu|^2 +u^2_\mu dx\leqslant\int_{\mathbb{R}^3}|\nabla \bar{u}|^2
 +\bar{u}^2
dx,
\end{equation}
and the lower semi-continuity of norm implies that
\begin{equation}\label{eq:5.7}
u_\mu\rightarrow \bar{u},\ \ \text{ in } H^1(\mathbb{R}^3) \ \text{
as }\mu\rightarrow+\infty.
\end{equation}
Now, we claim that $\bar{u}(x)\not\equiv0$ on $\Omega_0$.
Otherwise,  $\bar{u}=0$ a.e. in $x\in\mathbb{R}^3$, then
$|u_\mu|_{p+1}\overset{\mu\rightarrow+\infty}{\longrightarrow}0$
and
$\int\phi_{u_\mu}u_\mu^2dx\overset{\mu\rightarrow+\infty}{\longrightarrow}0$
by (\ref{eq:5.5.3.0}), hence $I(u_\mu)\rightarrow0$ as
$\mu\rightarrow+\infty$, which contradicts the fact that
$I(u_\mu)\geqslant\alpha>0$. Since $\partial\Omega_0$ is Lipschitz
continuous, we have $\bar{u}\in H_0^1(\Omega_0)$ by
$\bar{u}(x)=0$, a.e. in $\mathbb{R}^3\setminus\Omega_0$. It
follows from (\ref{eq:5.7}), (\ref{eq:5.5.3.1}) and (\ref{eq:5.5})
that, for any $\varphi\in H^1_0(\Omega_0)$,
\begin{equation}\nonumber
 \int_{\Omega_0}\nabla\bar{ u}\nabla\varphi+\bar{u}\varphi
 dx+\frac{\lambda}{4\pi}\int_{\Omega_0\times\Omega_0}\frac{\bar{u}^2(y)\bar{u}(x)\varphi(x)}{|x-y|}dydx
  =\int_{\Omega_0} |\bar{u}|^{p-1}\bar{u}\varphi dx,
\end{equation}
that is $\bar{u}$ is a weak solution of (\ref{eq:1.11}).
  $\Box$\\
{\bf Proof of Theorem \ref{th:1.5}:} This proof is motivated by that
of Theorem 1.3 in \cite{BartTPankovAWZhQ-CCM}. Let
$\mu_0=\max\{3(M_{0}^{p-1}-1),3(C_{p,\lambda}^{p-1}-1)\}$. By
Theorem \ref{th:1.1} or \ref{th:1.2}, the $L^\infty$-norms of
$u_\mu$ and $\phi_{u_\mu}$ are bounded uniformly in
$\mu\geqslant\mu_0$. By $\rm(G_3)$, there exists $R_0>0$ such that
$g(x)>\frac{5}{6}$ for $|x|\geqslant R_0$. Then, we have
\begin{equation}\label{eq:5.1.1}
(-\Delta+\mu)u_\mu^2=-2(W_\mu-\frac{\mu}{2})u_\mu^2-2|\nabla
u_\mu|^2\leqslant0\text{ for } |x|>R_0\text{ and }\mu>\mu_0,
\end{equation}
where $W_\mu(x)=1+\mu g(x)+\lambda\phi_{u_\mu}(x)-|u_\mu(x)|^{p-1}$.
On the other hand, let $w_\mu$ be the fundamental solution of
$-\Delta+\mu$, then, it follows from Proposition 3.1 and Theorem 4.2
in \cite{BerezinFA-ShubinMA} that $w_\mu$ is $C^\infty$ outside the
origin and it is positive such that
\begin{equation}\label{eq:5.1.2}
(-\Delta+\mu)w_\mu=0\text{ for } |x|>R_0\ \text{ and }
\end{equation}
\begin{equation}\label{eq:5.1.2.0}
w_\mu=|x|^{\frac{-(N-1)}{2}}e^{-\sqrt{\mu}|x|}(1+o(1)) \text{ as }
|x|\rightarrow+\infty, \text{ uniformly in } \mu\geqslant \mu_0.
\end{equation}
Let $C_0=\max\{M_0, C_{p,\lambda}\}$,
$\tilde{A}\geqslant\frac{C_0^2}{w_\mu(R_0)}$ and
$v(x)=u^2_\mu(x)-\tilde{A}w_\mu$. For $\mu>\mu_0$, it follows from
(\ref{eq:5.1.1}) and (\ref{eq:5.1.2}) that
\begin{equation}\nonumber
\left\{\begin{array}{ll}
 (-\Delta +\mu)v \leqslant0,\ \ \ |x|>R_0, \\
 v\leqslant0,\ \ \ \ |x|=R_0.
 \end{array}\right.
 \end{equation}
The maximum principle (Theorem 8.1 in \cite{DGilbargNStrudinger})
implies that
\begin{equation}\nonumber
u^2_\mu(x)\leqslant \tilde{A} w_\mu(x)\ \text{ for } |x|>R_0.
\end{equation}
By (\ref{eq:5.1.2.0}), there exists $A>0$ independent of $\mu>\mu_0$
such that
\begin{equation}\nonumber
u_\mu(x)\leqslant A
|x|^{-\frac{1}{2}}e^{-\frac{\sqrt{\mu}}{2}(|x|-R_0)}\ \text{ for }
|x|>R_0 \text{ and }\mu>\mu_0.\ \ \ \ \Box
\end{equation}

\section{Ground state of (\ref{eq:1.1}) for $p\in(1,2)$}

\begin{proposition}\label{Le:6.2} Let $\mathscr{N}=\{u\in \mathscr{D}_V\setminus\{0\}: I'(u)=0 \text{ in }
\mathscr{D}_V^\ast\}$ and $\lambda_\ast$ be given by Theorem
\ref{th:1.1}, $C_{p,\lambda}$ be given by Lemma \ref{Le:3.1}. Then,
for each $\lambda\in(0,\lambda_\ast)$ and
$\mu>C_{p,\lambda}^{p-1}-1$,
\begin{description}
 \item $\bf (i)$  $\mathscr{N}\neq\emptyset$, and there exists $C_p>0$ such that $\|u\|_{\mathscr{D}_V}\geqslant C_p$ for all $u\in \mathscr{N}$.
\item $\bf (ii)$ There exists $C_{\lambda,\mu}>0$ such that $\sup\limits_{u\in \mathscr{N}}(\|u\|_{\mathscr{D}_V} +|\nabla\phi_u|_2)\leqslant
C_{\lambda,\mu}$.
\item $\bf (iii)$  $\mathscr{N}\subset \mathscr{D}_V$ is a
compact set.
\end{description}
\end{proposition}
{\it Proof.} $\bf (i)$. Since $\lambda\in(0,\lambda^\ast)$ and
$\mu>C_{p,\lambda}^{p-1}-1$, Theorem \ref{th:1.1} implies that
$\mathscr{N}\neq\emptyset$. If $u\in\mathscr{N}$, then $I'(u)=0$ and
$u\not\equiv0$, that is,
\begin{equation}\nonumber
 \int_{\mathbb{R}^3}|\nabla u|^2+V_\mu
 u^2dx+\lambda\int_{\mathbb{R}^3}\phi_u
 u^2dx=\int_{\mathbb{R}^3}|u|^{p+1}dx.
\end{equation}
Hence,
$\|u\|_{\mathscr{D}_V}^2\leqslant\int_{\mathbb{R}^3}|u|^{p+1}dx\leqslant
S^{p+1} \|u\|_{\mathscr{D}_V}^{p+1}$ and
$\|u\|_{\mathscr{D}_V}\geqslant
S^{-\frac{p+1}{p-1}}:=C_p$.\\
$\bf (ii)$. Since $\|\cdot\|_{\mathscr{D}_V}$ and
$\|\cdot\|_{H^1(\mathbb{R}^3)}$ are equivalent norms in
$\mathscr{D}_V$, by Lemma \ref{Le:3.1} with $\Omega=\mathbb{R}^3$,
we have that $\sup\limits_{u\in \mathscr{N}}|u|_\infty \leqslant
C_{p,\lambda}$ and $|u(x)|\leqslant c_p\phi_u(x)$ a.e. in
$x\in\mathbb{R}^3$. For $u\in\mathscr{N}$, $I'(u)=0$ and
(\ref{eq:4.0.0}) holds.
 Then Lemma
\ref{Le:4.1} implies that
$$\|u\|_{\mathscr{D}_V}+|\nabla
\phi_u|_2=\|u\|_{\mathscr{D}_V}+\left\{\int_{\mathbb{R}^3}\phi_uu^2dx\right\}^{\frac{1}{2}}
<C_{\lambda,\mu}, \text{ for some } C_{\lambda,\mu}>0.$$
 $\bf (iii)$. Let $\{u_n\}\subset \mathscr{N}$, then $I'(u_n)=0$ and (\ref{eq:4.8}) holds for $u_n$.
By part $\bf (ii)$ and Lemma \ref{Le:4.2}, there exists
$u\in\mathscr{D}_V$ such that $u_n\overset{n}{\rightarrow}u$
 strongly in $L^{q}(\mathbb{R}^3)$. This and Lemma \ref{Le:2.2} show that
$\phi_{u_n}\overset{n}{\rightarrow}\phi_u$ strongly in
$\mathscr{D}^{1,2}(\mathbb{R}^3)$. Thus, for any $\varphi\in
\mathscr{D}_V$,
$$\int_{\mathbb{R}^3}\nabla u_n\nabla\varphi+V_\mu u_n\varphi dx
\overset{n}{\rightarrow}\int_{\mathbb{R}^3}\nabla
u\nabla\varphi+V_\mu u\varphi dx,$$
$$\int_{\mathbb{R}^3}\phi_{u_n}u_n\varphi dx
\overset{n}{\rightarrow}\int_{\mathbb{R}^3}\phi_{u}u\varphi dx\text{
and }\int_{\mathbb{R}^3}|u_n|^{p-1}u_n\varphi dx
\overset{n}{\rightarrow}\int_{\mathbb{R}^3}|u|^{p-1}u\varphi dx.$$
It follows from $I'(u_n)=0$ that $I'(u)=0$ in $\mathscr{D}_V^\ast$.
Then, using $I'(u_n)u_n=I'(u)u$ we know that
$u_n\overset{n}{\rightarrow}u$
 strongly in $\mathscr{D}_V$,
 this and $\bf (i)$ show that $\|u\|_{\mathscr{D}_V}\geqslant C_p$,
hence $u\in \mathscr{N}$ and $\mathscr{N}$ is a compact set in
$\mathscr{D}_V$. $\Box$\

{\bf Proof of Theorem \ref{th:6.1}.} By Proposition \ref{Le:6.2}
$\bf(iii)$, we know that $I$ is bounded below in $\mathscr{N}$ and
$c_0\triangleq\inf\limits_{u\in \mathscr{N}}I(u)> -\infty$, then
there exists $\{u_n\}\subset \mathscr{N}$ such that
$I(u_n)\overset{n}{\rightarrow}c_0$. So, Proposition \ref{Le:6.2}
$\bf(ii),(iii)$ imply that, there exists $u_0\in \mathscr{N}$ such
that $I(u_0)=c_0$ and  $I'(u_0)=0$, that is, $u_0$ is a ground state of (\ref{eq:1.1}).\ \ $\Box$ \\
\section{Appendix}
{\footnotesize
 {\bf Proof of Lemma \ref{Le:3.2}:} For $\beta>1$, $M>1$, the same
as \cite{TNS-ASNSPisa-1986}, we define a function $H\in
C^1(0,+\infty)$ by
$$H(s)=\left\{\begin{array}{ll}
 s^\beta,  \hfill s\in [0,M], \\
 \beta M^{\beta-1}s-(\beta-1)M^{\beta},\  \hfill s\in
 (M,\infty).
 \end{array}\right.$$
 Let
$$G(s)\overset{\vartriangle}{=}\int_{0}^{s}|H{'}(t)|^{2}dt=\left\{\begin{array}{ll}
 \frac{\beta^2}{2\beta-1}s^{2\beta-1},  \hfill s\in (0,M], \\
 \beta^2M^{2(\beta-1)}s-\frac{2\beta^2(\beta-1)}{2\beta-1}M^{2\beta-1},  \hfill s\in
 (M,\infty).
 \end{array}\right.$$
 Then $G(s)$ and $H(s)$ are Lipschitz in $[0,+\infty)$, therefore,
 $G(u)$ and $H(u)$ are in $\mathscr{D}^{1,2}(\mathbb{R}^N)$ if $u\in
\mathscr{D}^{1,2}(\mathbb{R}^N)$. Moreover,
 \begin{equation}\label{eq:3.6}
sG(s)\leqslant s^2 {H{'}}^2(s)\leqslant\beta^{2}H^{2}(s).
\end{equation}
Let $\bar{\eta}(x)\in C^\infty_0(\mathbb{R}^N)$ with
$\bar{\eta}(x)\geqslant0$ such that, for $r_1>r_2>0$( $r_1$ will be
determined later)
$$\bar{\eta}(x)\equiv1\text{ for }x\in B_{r_2},\ \ \bar{\eta}(x)\equiv0\ \text{ for }x\in B^c_{r_1},
 \ \text{and }|\nabla\bar{\eta}|\leqslant\frac{2}{r_1-r_2}.$$
for each $y\in\mathbb{R}^N$, setting
$\eta(x)=\bar{\eta}(x-y)\geqslant0$ and $0\leqslant\eta^2G(u)\in
\mathscr{D}^{1,2}(\mathbb{R}^N)$ with compact support in
$B_{r_1}(y)\cap\{x\in\mathbb{R}^3:u(x)\neq0\}$. By (\ref{eq:3.6})
and  (\ref{eq:3.5}) with $\varphi=\eta^2$ and $h=G$,
\begin{eqnarray}\label{eq:3.7}
\int_{\mathbb{R}^{N}}\nabla u\nabla (\eta^2 G(u))
dx\leqslant\int_{\mathbb{R}^{N}}u^{p-1}u\eta^2 G(u)
dx\leqslant\beta^2\int_{\mathbb{R}^{N}}|u|^{p-1}\eta^2 H^2(u) dx.
\end{eqnarray}
Noting that
\begin{equation}\nonumber
\nabla u\nabla(\eta^2 G(u))=|\nabla u|^2{H{'}}^2(u)\eta^2+2\nabla
u\nabla\eta G(u)\eta.
\end{equation}
Also, by Young inequality and (\ref{eq:3.6}),
\begin{equation}\nonumber
\begin{split}
\left|\nabla u\nabla\eta G(u)\eta\right|&=|\eta
u^{-1/2}G^{1/2}(u)\nabla u||
u^{1/2}G^{1/2}(u)\nabla\eta|\\
&\leqslant\frac{1}{4}|\nabla
u|^2{H{'}}^2(u)\eta^2+4\beta^2{|\nabla\eta|}^2 H^2(u).
\end{split}
\end{equation}
Then, it follows from (\ref{eq:3.7}) that
\begin{eqnarray}\nonumber
\int_{\mathbb{R}^{N}}|\nabla u|^2|H{'}(u)|^2\eta^2 dx\leqslant
16\beta^2\left(\int_{\mathbb{R}^{N}}{|\nabla\eta|}^2 {H^2(u)}
dx+\int_{\mathbb{R}^{N}}|u|^{p-1}\eta^2 H^2(u)  dx\right).
\end{eqnarray}
Hence, by $\beta>1$ and H\"{o}lder inequality we see that
\begin{equation}\label{eq:3.8}
\begin{split}
\int_{\mathbb{R}^{N}}|\nabla({H(u)}\eta)|^2 dx
&\leqslant2\int_{\mathbb{R}^N}|\nabla\eta|^2H^2(u)dx+2\int_{\mathbb{R}^N}|\nabla u|^2|H'(u)|^2\eta^2dx \\
&\leqslant 34\beta^2\left(\int_{\mathbb{R}^{N}}{|\nabla\eta|}^2
H^2(u) dx+\int_{\mathbb{R}^{N}}|u|^{p-1}\eta^2 H^2(u)  dx\right)\\
&\leqslant 34\beta^2\left(\int_{\mathbb{R}^{N}}{|\nabla\eta|}^2
H^2(u) dx+|\eta
H(u)|_{2^\ast}^2|u|^{p-1}_{2^\ast}|B_{r_1}|^{\gamma}\right),
\end{split}
\end{equation}
where $\gamma=\frac{2^\ast-p-1}{2^\ast}$ and $|B_{r_1}|$ denotes the
volume of $B_{r_1}$. Taking
$\beta=\frac{2^\ast}{2}\triangleq\beta_0$ and
\begin{equation}\label{eq:3.9}
r_1=|B_1|^{-\frac{1}{N}}\left(68\beta_0^2|u|_{2^\ast}^{p-1}+1\right)^{-\frac{\gamma}{N}}
\leqslant\min\left\{|B_1|^{-\frac{1}{N}},|B_1|^{-\frac{1}{N}}\left(68\beta_0^2|u|_{2^\ast}^{p-1}\right)^{-\frac{\gamma}{N}}\right\}.
\end{equation}
Obviously,
$34\beta_0^2|u|^{p-1}_{2^\ast}|B_{r_1}|^{\gamma}\leqslant\frac{1}{2}$
and $|B_{r_1}|\leqslant1$. Then, (\ref{eq:3.8}) with $\beta=\beta_0$
gives that
\begin{equation}\nonumber
\begin{split}
\int_{\mathbb{R}^{N}}|\nabla({H(u)}\eta)|^2 dx \leqslant
68\beta_0^2\int_{\mathbb{R}^{N}}{|\nabla\eta|}^2 H^2(u) dx.
\end{split}
\end{equation}
Hence, by Sobolev inequality and $\eta H(u)\in
\mathscr{D}_0^{1,2}(B_{r_1}(y))$ as well as the definition of
$\eta$, we see that
\begin{eqnarray}\label{eq:3.10}
|H(u)|^{2}_{L^{2^\ast}(B_{r_2}(y))}&\leqslant&\left(\int_{\mathbb{R}^{N}}|H(u)\eta|^{2^\ast}dx\right)^{2/2^\ast}
\leqslant
S_0^{-1}\int_{\mathbb{R}^{N}}|\nabla({H(u)}\eta)|^2 dx\nonumber\\
&\leqslant& 68S_0^{-1}\beta_0^2\int_{\mathbb{R}^{N}}|\nabla\eta|^2
H^2(u)dx\nonumber\\
&\leqslant& \left(\frac{C\beta_0}{r_1-r_2}\right)^2|
H(u)|^2_{L^{2}(B_{r_1}(y))},
\end{eqnarray}
where $C=\sqrt{68S_0^{-1}}$ and $S_0$ is the Sobolev constant, which
depends only on $N$. By definition of $H(s)\equiv s^\beta$ if
$M\rightarrow+\infty$. Noting that $\beta=\beta_0$ and
$2\beta_0=2^\ast$, it follows from (\ref{eq:3.10}) that
\begin{eqnarray}\label{eq:3.11}
|u|^{2\beta_0}_{L^{2\beta_0^2}(B_{r_2}(y))}\leqslant
\left(\frac{C\beta_0}{r_1-r_2}\right)^2|u|^{2\beta_0}_{L^{2^\ast}(B_{r_1}(y))}.
\end{eqnarray}
For each $i\geqslant2$, let $r_i=\frac{2+2^{-i}}{4}r_1$ ( $r_1$ is
given by (\ref{eq:3.9})). Take $\bar{\eta}_i\in
C^\infty_0(\mathbb{R}^3)$ such that
$$\bar{\eta}_i(x)\equiv1\text{ for }x\in B_{r_{i+1}},\ \ \bar{\eta}_i(x)\equiv0\ \text{ for }x\in \mathbb{R}^N\setminus B_{r_i}
 \ \text{and }|\nabla\bar{\eta}_i|\leqslant\frac{2}{r_i-r_{i+1}}.$$
Let $\delta=\frac{2\beta_0}{2\beta_0^2+1-p}$ and $\delta\in(0,1)$ by
$p\in(1,\frac{N+2}{N-2})$. For $i\geqslant2$, applying
(\ref{eq:3.8}) with $\eta=\eta_i(x)\triangleq\bar{\eta}_i(x-y)$ and
$\beta=\beta_i\triangleq\delta^{-i}>1$, then noting that
$|B_{r_i}|<|B_{r_1}|\leqslant1$ and using H\"{o}lder inequality and
the definition of $\bar{\eta}_i$, we get
\begin{equation}\nonumber
\begin{split}
&\int_{\mathbb{R}^{N}}|\nabla({H(u)}\eta_i)|^2 dx \leqslant
34\beta_i^2\left(\int_{\mathbb{R}^{N}}{|\nabla\eta_i|}^2
H^2(u) dx+\int_{\mathbb{R}^{N}}|u|^{p-1}\eta_i^2 H^2(u)  dx\right)\\
&\leqslant
34\beta_i^2\left[\left(\frac{2}{r_i-r_{i+1}}\right)^2|H(u)|^{2}_{L^{2^{\ast}\delta}(B_{r_i}(y))}
+|H(u)|^{2}_{L^{2^{\ast}\delta}(B_{r_i}(y))}|u|^{p-1}_{L^{2\beta_0^{2}}(B_{r_i}(y))}\right]\\
&\leqslant 34\beta_i^2\left[\left(\frac{2}{r_i-r_{i+1}}\right)^2
+|u|^{p-1}_{L^{2\beta_0^{2}}(B_{r_2}(y))}\right]|H(u)|^{2}_{L^{2^{\ast}\delta}(B_{r_i}(y))}\\
 &\leqslant
34\beta_i^2\left[\left(\frac{2}{r_i-r_{i+1}}\right)^2
+\left(\frac{2}{r_1-r_{2}}\right)^{(p-1)/\beta_0}|u|^{p-1}_{L^{2^\ast}(B_{r_1}(y))}\right]|H(u)|^{2}_{L^{2^{\ast}\delta}(B_{r_i}(y))}
\text{ by (\ref{eq:3.11}) }\\
 &\leqslant
\left(\frac{\bar{C}\beta_i(1+|u|_{2^\ast}^{(p-1)/2})}{r_i-r_{i+1}}\right)^2|H(u)|^{2}_{L^{2^{\ast}\delta}(B_{r_i}(y))},
\end{split}
\end{equation}
where $\bar{C}$ is a constant depending only on $N$ and $p$. Since
$\eta_i H(u)\in\mathscr{D}^{1,2}_0(B_{r_i}(y))$, it follows from the
Sobolev inequality that
\begin{eqnarray}\nonumber
|H(u)|^{2}_{L^{2^\ast}(B_{r_{i+1}}(y))} \leqslant
S_0^{-1}\int_{\mathbb{R}^{N}}|\nabla({H(u)}\eta_i)|^2 dx\leqslant
\left(\frac{\bar{C}_1\beta_i(1+|u|_{2^\ast}^{(p-1)/2})}{r_i-r_{i+1}}\right)^2|H(u)|^{2}_{L^{2^{\ast}\delta}(B_{r_i}(y))}
\end{eqnarray}
where $\bar{C}_1=\bar{C}\sqrt{S_0^{-1}}$ depends only on $N$ and
$p$. Let $M\rightarrow+\infty$, and $H(u)=u^{\beta_i}$, then
\begin{eqnarray}\nonumber
|u|_{L^{2^\ast\beta_i}(B_{r_{i+1}}(y))}\leqslant
\left(\frac{\bar{C}_1\beta_i(1+|u|_{2^\ast}^{(p-1)/2})}{r_i-r_{i+1}}\right)^{1/\beta_i}|u|_{L^{2^{\ast}\beta_{i-1}}(B_{r_i}(y))}.
\end{eqnarray}
We can now perform the Moser iterations in a standard way and get
that
\begin{eqnarray}\label{eq:3.12}
|u|_{L^{2^\ast\beta_i}(B_{r_{i+1}}(y))}&\leqslant&
\prod_{l=2}^{i}\left(\frac{\bar{C}_1\beta_l(1+|u|_{2^\ast}^{(p-1)/2})}{r_l-r_{l+1}}\right)^{\frac{1}{\beta_l}}|u|_{L^{2^\ast\beta_1}(B_{r_2}(y))}\\
 &=&(\frac{2}{\delta})^{f(i)}\left(\frac{8\bar{C}_1(1+|u|_{2^\ast}^{(p-1)/2})}{r_1}\right)^{g(i)}\!\!|u|_{L^{2^\ast\beta_1}(B_{r_2}(y))}\nonumber.
\end{eqnarray}
where
$f(i)=\frac{2\delta^2}{1-\delta}+\frac{\delta^3(1-\delta^{i-2})}{(1-\delta)^2}
+\frac{i\delta^{i+1}}{1-\delta}\overset{i}{\rightarrow}\frac{2\delta^2-\delta^3}{(1-\delta)^2}$
and
$g(j)=\frac{\delta^2(1-\delta^{i-1})}{1-\delta}\overset{i}{\rightarrow}\frac{\delta^2}{1-\delta}$.
Since $2^\ast\beta_1=2\beta_0^2+1-p<2\beta_0^2$, by H\"{o}lder
inequality, (\ref{eq:3.9}) and (\ref{eq:3.11}), we have
\begin{eqnarray}\label{eq3.12}
|u|_{L^{2^\ast\beta_1}(B_{r_2}(y))}
\leqslant|B_{r_2}|^{\frac{2\beta_0^2}{p-1}}|u|_{L^{2\beta_0^2}(B_{r_2}(y))}\leqslant|u|_{L^{2\beta_0^2}(B_{r_2}(y))}
\leqslant\left(\frac{C\beta_0}{r_1-r_2}\right)^{1/\beta_0}|u|_{2^\ast}\nonumber.
\end{eqnarray}
This and (\ref{eq:3.12}) show that
\begin{eqnarray}\nonumber
|u|_{L^{2^\ast\beta_i}(B_{r_1/2}(y))}\leqslant|u|_{L^{2^\ast\beta_i}(B_{r_{i+1}}(y))}
 \leqslant(\frac{2}{\delta})^{f(i)}\left(\frac{8\bar{C}_1(1+|u|_{2^\ast}^{(p-1)/2})}{r_1}\right)^{g(i)}
 \left(\frac{\frac{16}{7}C\beta_0}{r_1}\right)^{1/\beta_0}|u|_{2^\ast}\nonumber.
\end{eqnarray}
Let $i\rightarrow+\infty$, and noting (\ref{eq:3.9}), There exist
positive constant $C_1(p,N),C_2(p,N)$ such that
\begin{eqnarray}\nonumber
 &&|u|_{L^{\infty}(B_{r_1/2}(y))}
\leqslant(\frac{2}{\delta})^{\frac{2\delta^2-\delta^3}{(1-\delta)^2}}\left(\frac{8\bar{C}_1(1+|u|_{2^\ast}^{(p-1)/2})}{r_1}\right)^{\frac{\delta^2}{1-\delta}}
 \left(\frac{\frac{16}{7}C\beta_0}{r_1}\right)^{1/\beta_0}|u|_{2^\ast}\nonumber\\
 &\leqslant&(\frac{2}{\delta})^{\frac{2\delta^2-\delta^3}{(1-\delta)^2}}\left(\frac{(8\bar{C}_1+\frac{16}{7}C\beta_0)
 (1+|u|_{2^\ast}^{\frac{p-1}{2}})}{r_1}\right)^{\frac{\delta^2}{1-\delta}+\frac{1}{\beta_0}}
 |u|_{2^\ast}\nonumber\\
 &\leqslant&(\frac{2}{\delta})^{\frac{2\delta^2-\delta^3}{(1-\delta)^2}}
 \left(|B_1|^{\frac{1}{N}}(8\bar{C}_1+\frac{16}{7}C\beta_0)
 (1+|u|_{2^\ast}^{\frac{p-1}{2}})
 \left(68\beta_0^2|u|_{2^\ast}^{p-1}+1\right)^{\frac{\gamma}{N}}\right)^{\frac{\delta^2}{1-\delta}+\frac{1}{\beta_0}}
 |u|_{2^\ast}\nonumber\\
 &\leqslant&C_1(p,N)
 (1+|u|_{2^\ast}^{C_2(p,N)})
 |u|_{2^\ast}\nonumber.
\end{eqnarray}
Since $y\in\mathbb{R}^N$ and $r_1$ is fixed by (\ref{eq:3.9}), we
have
\begin{eqnarray}\nonumber
|u|_{\infty}\leqslant
C_1(p,N)\left(1+|u|_{2^\ast}^{C_2(p,N)}\right)|u|_{2^\ast}.\ \ \
\Box
\end{eqnarray}}

\end{document}